\documentclass{article}

\usepackage{amsmath}
\usepackage{amssymb}
\usepackage{amsthm,mathtools}
\usepackage{bbm}
\usepackage{hyperref}
\usepackage{theoremref}
\usepackage{enumitem}

\newcommand{\ap}[1]{\prescript{a}{}{#1}}

\def\reals{\mathbbm{R}}
\def\ereals{\overline{\reals}}

\def\uball{\mathbbm{B}}

\def\comp{\mathop{\text{\scriptsize\raise 1pt \hbox{$\circ$}}}}
\def\infconv{\mathop{\text{\scriptsize\raise 1pt \hbox{$\square$}}}}

\def\minimize{\mathop{\rm minimize}\limits}
\def\maximize{\mathop{\rm maximize}\limits}

\def\st{\mathop{\rm subject\ to}}

\def\dom{\mathop{\rm dom}\nolimits}
\def\gph{\mathop{\rm gph}\nolimits}	
\def\core{\mathop{\rm core}\nolimits}
\def\rcore{\mathop{\rm rcore}\nolimits}

\def\lev{\mathop{\rm lev}\nolimits}
\def\ovr{\mathop{\rm over}\ }
\def\rge{\mathop{\rm rge}}

\def\upto{{\raise 1pt \hbox{$\scriptstyle \,\nearrow\,$}}}
\def\downto{{\raise 1pt \hbox{$\scriptstyle \,\searrow\,$}}}
\def\inte{\mathop{\rm int}}
\def\rinte{\mathop{\rm rint}\nolimits}

\def\pos{\mathop{\rm pos}}
\def\aff{\mathop{\rm aff}\nolimits}

\def\epi{\mathop{\rm epi}\nolimits}
\DeclareMathOperator{\sign}{sign}

\def\tos{\rightrightarrows}
\def\FF{(\F_t)_{t=0}^T}
\def\one{\mathbbm 1}

\def\A{{\cal A}}

\def\C{{\cal C}}
\def\D{{\cal D}}

\def\F{{\cal F}}
\def\G{{\cal G}}

\def\K{{\cal K}}
\def\L{{\cal L}}
\def\M{{\cal M}}
\def\N{{\cal N}}

\def\S{{\cal S}}

\def\T{{\cal T}}
\def\U{{\cal U}}
\def\V{{\cal V}}
\def\X{{\cal X}}
\def\Y{{\cal Y}}

\newtheorem{theorem}{Theorem}
\newtheorem{lemma}[theorem]{Lemma}
\newtheorem{corollary}[theorem]{Corollary}

\newtheorem{example}[theorem]{Example}
\newtheorem{cexample}[theorem]{Counterexample}
\newtheorem{remark}[theorem]{Remark}

\theoremstyle{definition}

\newtheorem{assumption}[theorem]{Assumption}

\title{Dual solutions in convex stochastic optimization}

\author{Teemu Pennanen\thanks{Department of Mathematics, King's College London, Strand, London, WC2R 2LS, United Kingdom, teemu.pennanen@kcl.ac.uk} \and Ari-Pekka Perkki\"o\thanks{Mathematics Institute, Ludwig-Maximilian University of Munich, Theresienstr. 39, 80333 Munich, Germany, a.perkkioe@lmu.de. Corresponding author}}

\begin{document}

\maketitle

\begin{abstract}
This paper studies duality and optimality conditions for general convex stochastic optimization problems. The main result gives sufficient conditions for the absence of a duality gap and the existence of dual solutions in a locally convex space of random variables. It implies, in particular, the necessity of scenario-wise optimality conditions that are behind many fundamental results in operations research, stochastic optimal control and financial mathematics. Our analysis builds on the theory of Fr\'echet spaces of random variables whose topological dual can be identified with the direct sum of another space of random variables and a space of singular functionals. The results are illustrated by deriving sufficient and necessary optimality conditions for several more specific problem classes. We obtain significant extensions to earlier models e.g.\ on stochastic optimal control, portfolio optimization and mathematical programming.
\end{abstract}


\section{Introduction}

Given a probability space $(\Omega,\F,P)$ with a filtration $\FF$ (an increasing sequence of sub-$\sigma$-algebras of $\F$), consider the problem
\begin{equation}\label{sp}\tag{$SP$}
  \begin{aligned}
   &\minimize\quad & & Ef(x,\bar u)\quad\ovr\text{$x\in\N$}
  \end{aligned}
\end{equation}
where $\N$ is a linear space of stochastic processes $x=(x_t)_{t=0}^T$ adapted to $\FF$ (i.e., $x_t$ is $\F_t$-measurable) and $\bar u$ is a $\reals^m$-valued random variable. We assume that $x_t$ takes values in a Euclidean space $\reals^{n_t}$ so the process $x=(x_t)_{t=0}^T$ takes values in $\reals^n$ where $n:=n_0+\cdots+n_T$. The objective is defined on the space $L^0(\Omega,\F,P;\reals^n\times\reals^m)$ of $\F$-measurable $\reals^n\times\reals^m$-valued random variables $(x,u)$ by
\[
Ef(x,u) := \int_\Omega f(x(\omega),u(\omega),\omega)dP(\omega),
\]
where $f$ is a {\em convex normal integrand} on $\reals^n\times\reals^m\times\Omega$, i.e.\ $f(\cdot,\omega)$ is a closed convex function for every $\omega\in\Omega$ and $\omega\mapsto\epi f(\cdot,\omega)$ is an $\F$-measurable set-valued mapping; see \cite[Chapter~14]{rw98}. Here and in what follows, we define the integral of an extended real-valued random variable as $+\infty$ unless its positive part is integrable. By \cite[Proposition~14.28]{rw98}, the function $\omega\mapsto f(x(\omega),\bar u(\omega),\omega)$ is measurable for any measurable $x$ and $\bar u$ so the integral functional $Ef$ is a well-defined convex function on $L^0(\reals^n\times\reals^m)$.

Problems of the form \eqref{sp} unify and extend various more specific models in operations research, engineering and economics and they have been extensively studied since they first appeared in \cite{roc76}; see e.g.\ \cite{rw78,rw83,pli82,sdr9,pen11c,cccd15,bpp18} and the references there. Like \cite{rw76,pli82,pen11c,bpp18}, this paper studies \eqref{sp} through convex duality. We follow the recent approach of \cite{pp22b}, which yields explicit dual problems and allows for general adapted strategies $x$ without any integrability restrictions. This is important e.g.\ in common models of financial mathematics where it may be impossible to find optimal trading strategies in an integrable space of stochastic processes; see \cite{pp22b} and its references.

The analysis of \cite{pp22b} is based on applying the conjugate duality framework of \cite{roc74} to the parametric optimization problem
\begin{equation}\label{spx}\tag{$SP_\X$}
\begin{aligned}
  &\minimize\quad Ef(x,\bar u)\quad\ovr x\in\X\\
  &\st\quad x-z\in\X_a
\end{aligned}
\end{equation}
where $\X\subset L^0(\Omega,\F,P:\reals^n)$ is a locally convex vector space of random variables,  
\[
\X_a:=\X\cap\N
\]
and $z$ and $\bar u$ are the parameters. We assume that $z\in\X$ and that $\bar u$ belongs to another locally convex space $\U\subset L^0(\Omega,\F,P;\reals^m)$ of random variables. Following \cite{pp22b}, we apply the conjugate duality framework to the extended real-valued convex function $F$ on $\X\times\X\times\U$ defined by
\[
F(x,z,u) := Ef(x,u) + \delta_\N(x-z).
\]
The associated optimum value function is denoted by
\[
\varphi(z,u):=\inf_{x\in\X}\{Ef(x,u)\,|\,x-z\in\N\}.
\]
We assume that $\X$ and $\U$ are in separating duality with spaces of random variables $\V\subset L^0(\Omega,\F,P;\reals^n)$ and $\Y\subset L^0(\Omega,\F,P;\reals^m)$, respectively. The bilinear forms are the usual ones, i.e.
\[
\langle x,v\rangle:=E[x\cdot v]\quad\text{and}\quad \langle u,y\rangle := E[u\cdot y].
\]
It is assumed that all the spaces are {\em decomposable} and {\em solid}. Decomposability of $\X$ means that
\[
1_Ax+1_{\Omega\setminus A}x'\in\X
\]
for every $x\in\X$ and $x'\in L^\infty$ while solidity means that if $\bar x\in\X$ and $x\in L^0$ are such that $|x^i|\le|\bar x^i|$ almost surely for every $i=1,\ldots,m$, then $x\in\X$. Similarly for $\U$, $\V$ and $\Y$.

As soon as $\dom Ef\cap(\X\times\U)\ne\emptyset$, the conjugate of $F$ can be expressed as
\[
F^*(v,p,y) = Ef^*(v+p,y)+\delta_{\X_a^\perp}(p);
\]
see \cite[Lemma~13]{pp22b}. It follows that the dual problem becomes
\begin{equation}\label{d}\tag{$D$}
\maximize\quad \langle\bar u,y\rangle - Ef^*(p,y)\quad\ovr\quad (p,y)\in\X_a^\perp\times\Y.
\end{equation}
It was shown in \cite{pp22b} that if the optimum value function $\varphi$ is subdifferentiable at $(0,\bar u)$ then there is no duality gap and a feasible $x\in\N$ solves \eqref{sp} if and only if there exists a dual feasible $(p,y)\in\X_a^\perp\times\Y$ such that the scenario-wise optimality condition
\[
(p,y)\in\partial f(x,\bar u)\text{ a.s.}
\]
holds. It should be noted that the above applies to the original problem \eqref{sp} which does not directly fit the conjugate duality framework of \cite{roc74} as the space $\N$ is not locally convex. Allowing for general adapted strategies is key to obtaining the existence of primal solutions and the lower semicontinuity of the optimum value function $\varphi$ in applications where the feasible set may be unbounded; see \cite[Section~5]{pp22b}. An example is the celebrated ``fundamental theorem of asset pricing'' in financial mathematics \cite{dmw90,sch92}.

This paper gives sufficient conditions for subdifferentiability of the optimum value function $\varphi$ of \eqref{spx}. A simple sufficient condition is that $\varphi$ be Mackey-continuous at $(0,\bar u)$ since that implies the existence of a subgradient in $\V\times\Y$; see \thref{relsd}. In some applications this is easy to establish e.g.\ by showing that $\varphi$ has a Mackey continuous upper bound; see \thref{cor:cont0} as well as the discussion after \thref{thm:fm5}. In others, however, the Mackey continuity fails. This is typical in problems with pointwise constraints such as those ones studied in \cite{pp22b} and in Section~\ref{sec:app} below.

We follow the main arguments from \cite{rw76,pp18a} and first establish, in Section~\ref{sec:sdp}, (relative) continuity of the value function $\varphi$ with respect to a Fr\'echet topology which may be strictly stronger than the Mackey topology. Strengthening the topology makes it possible to establish the continuity in many applications. However, continuity with respect to such a topology does not, in general, give subgradients in the Mackey-dual $\V\times\Y$ of $\X\times\U$ but in a strictly larger dual space. Under an additional condition that extends classical ``relatively complete recourse'' condition from \cite{rw76}, we will show, in Section~\ref{sec:recourse}, that the strong subdifferential of $\varphi$ at $(0,\bar u)$ intersects the space $\V\times\Y$ thus giving a subgradient with respect to the original dual pairing. As noted earlier, this implies the absence of a duality gap and the existence of solutions to the dual problem \eqref{d}.

Our main result extends those of \cite{rw76,pp18a} in two important says. First, we include the parameter $u\in\U$ which allows for dualization of different formats of stochastic optimization problems in a unified manner much like the conjugate duality framework of \cite{roc74} unifies more specific problem formats. Second, we go beyond spaces of essentially bounded strategies and parameters. Allowing for more general Fr\'echet spaces $\X$ and $\U$, widens the scope of applicability of the theory. The extensions are illustrated in Section~\ref{sec:app} by deriving optimality conditions in mathematical programming, stochastic optimal control and portfolio optimization in a generality not seen before. In particular, we extend the main results of \cite{rw78} by allowing equality constraints and unbounded feasible sets. In stochastic optimal control, we find necessary and sufficient optimality conditions under far more general conditions that have been given before. In particular, we allow for both state and control constraints and unbounded feasible sets which are encountered e.g.\ in the classical models of linear-quadratic control. In the case of portfolio optimization, we extend earlier results by the inclusion of portfolio constraints and statically held derivative portfolios. 

The strategy outlined above is based on endowing the spaces $\X$ and $\U$ with topologies under which they become Fr\'echet spaces and under which their topological duals can be expressed as direct sums of $\V$ and $\Y$, respectively, with linear spaces of certain singular functionals. Section~\ref{sec:ifst} reviews the basic theory of such spaces and extends earlier results on conjugates of integral functionals with respect to the associated dual pairings.

\section{Integral functionals under the strong topology}\label{sec:ifst}




This section studies convex integral functionals with respect to a ``strong topology'' on $\U$ which may be strictly stronger than the Mackey topology $\tau(\U,\Y)$. The stronger the topology, the easier it is to establish continuity of a functional. By \thref{relsd}, continuity implies subdifferentiability which, as we have seen, is closely connected to existence of dual solutions; see Section~\ref{sec:cs}. Particularly convenient situation is where $\U$ is a Fr\'echet space under the strong topology since then, lower semicontinuous convex functions on $\U$ are continuous throughout the {\em core} of their domain; see \thref{thm:cont0}.  Recall that a {\em Fr\'echet space} is a complete metrizable locally convex vector space. In particular, Banach spaces are Fr\'echet spaces.

Given a convex set $C\subset\U$, its core is the set
\[
\core C = \{u\in C\mid \pos(C-x)=\U\},
\]
where
\[
\pos(C-u):=\bigcup_{\alpha>0}\alpha(C-u)
\]
is the {\em positive hull} of $C-u$. More generally, \thref{thm:ipcont} says that if a convex function $g$ on a Fr\'echet space is the infimal projection of a closed convex function (as is the case with the optimum value function $\varphi$ of \eqref{sp}) and if the affine hull of $\dom g$ is closed, then $g$ is relatively continuous throughout the relative core of its domain. By \thref{relsd}, relative continuity is sufficient for subdifferentiability. Recall that the {\em relative core} of a set $C$ is its core with respect to its affine hull. By \thref{lem:rcore},
\[
\rcore C = \{u\in C\mid \pos(C-x) \text{ is linear}\}
\]
and, for every $u\in\core C$, the positive hull $\pos(C-u)$ is the linear translation of $\aff C$. In many applications, one encounters convex functions that fail to be continuous throughout the whole space but, nevertheless, are relatively continuous. 


If the strong topology is strictly stronger than the Mackey topology $\tau(\U,\Y)$, then the corresponding topological dual of $\U$ is strictly larger than $\Y$. The extra elements in the dual space cannot be represented by measurable functions but, often, they have properties that can be employed in the duality theory of integral functionals and, as we will see, in establishing the existence of solutions to the dual problem \eqref{d}.

\subsection{Continuity of integral functionals}\label{sec:cif}

We will continue with the assumption that $\U$ is a solid decomposable space of random variables. From now on, we will also assume that $\U$ is endowed with a topology $s$ at least as strong as the Mackey topology $\tau(\U,\Y)$ such that $s$ makes $\U$ Fr\'echet space. We will call $s$ the {\em strong topology}.

We say that a convex function $g$ is {\em relatively continuous} at a point $u\in\dom g$ if it is continuous at $u$ in the relative topology of $\aff\dom g$. The following is a direct consequence of \thref{cor:relcontiff}.

\begin{theorem}\thlabel{thm:Ehcont}
Assume that $Eh$ is proper and strongly lsc and that $\aff\dom Eh$ is strongly closed. Then $Eh$ is strongly relatively continuous throughout $\rcore\dom Eh$.
\end{theorem}

Recall that the  {\em K\"othe dual} or {\em associate space} of $\U$ is the linear space
\[
\U':=\{y\in L^0\mid u\cdot y\in L^1\quad \forall u\in\U\};
\]
see, e.g., \cite{bs88,zaa83}. We denote by $\U^*$ the topological dual of $\U$ under the strong topology. That is, $\U^*$ is the linear space of all strongly continuous linear functionals on $\U$.

\begin{lemma}\thlabel{lem:kothe}
We have $\U'\subseteq\U^*$ in the sense that, for every $y\in\U'$, the linear functional
\[
u\mapsto E[u\cdot y]
\]
is strongly continuous on $\U$.
\end{lemma}

\begin{proof}
Let $y\in L^0$ be such that $u\cdot y\in L^1$ for all $u\in\U$. Given $u\in\U$, the solidity of $\U$ implies that the random vector $u'\in L^0$ defined by $u'_i:=|u_i|\sign(y_i)$ belongs to $\U$. Since $u'\cdot y=\sum_{i=1}^m|u_i||y_i|$, the function
\[
p_y(u):=E[\sum_{i=1}^m|u_i||y_i|]
\]
is thus finite on $\U$. By Fatou's lemma, $p_y$ is $L^1$-lsc on $L^1$ so, by \cite[Lemma~6]{pp12}, it is $\sigma(\U,\Y)$-lsc. By \thref{thm:Ehcont}, $p_y$ is strongly continuous on $\U$. Since the function $u\mapsto E[u\cdot y]$ is majorized by $p_y$, it is strongly continuous by \thref{thm:cont0}.
\end{proof}



%


%

The following gives a simple criterion for strong continuity of an integral functional on an Orlicz space.

\begin{example}\thlabel{ex:orliczcont}
Let $\Phi:\reals\times\Omega\to\ereals$ be a convex normal integrand such that, for almost every $\omega\in\Omega$, the function $\Phi(\cdot,\omega)$ is nonconstant, symmetric, vanishes at the origin, has $\dom\Phi(\cdot,\omega)\ne\{0\}$ and $\Phi(a)\in L^1$ for some constant $a>0$. The Musielak-Orlicz space
\[
L^\Phi := \{x\in L^0\mid \exists \alpha>0: E\Phi(|x|/\alpha)<\infty\}
\]
is a solid decomposable Banach space; see \cite{pp220}. 

Let $h$ be a convex normal integrand such that $Eh$ is proper and lsc on $L^\Phi$ and there exist constants $\epsilon,M>0$ with
\[
h(\epsilon u)\le M\Phi(|u|)+M\quad\forall u\in\reals^n\ \text{a.s.}
\]
Then  $0\in\core\dom Eh$ and $Eh$ is continuous throughout $\core\dom Eh$. Defining
\[
\tilde h(u,\omega):=h(u,\omega)+\delta_{L(\omega)}(u),
\]
where $L$ is a random linear set, we have $0\in\rcore\dom E\tilde h$,
\[
\aff\dom E\tilde h = L^\Phi(L) = L^\Phi(\aff\dom\tilde h)
\]
and $E\tilde h$ is strongly relatively continuous throughout $\rcore\dom E\tilde h$.
\end{example}

\begin{proof}
The continuity on $\core\dom Eh$ follows from \thref{thm:Ehcont}. Given $u\in L^\Phi$, there exists an $\alpha>0$ such that $E\Phi(|u|/\alpha)\le 1$ so the assumption gives the existence of an $\epsilon>0$ such that $Eh(\epsilon u)\le 2M$. Thus, $0\in\core\dom Eh$.

The claims concerning $\tilde h$ follow from the fact that $E\tilde h=Eh+\delta_{L^\Phi(L)}$, where
\[
L^\Phi(L):=\{u\in L^\Phi\mid u\in L\ a.s.\}
\]
is strongly closed, since it is $\sigma(\U,\Y)$-closed, by \cite[Corollary~7]{pp22b}. Indeed, by \thref{lem:hulls},
\[
\pos(\dom E\tilde h) = \pos(\dom E h\cap L^\Phi(L)) = \pos(\dom E h)\cap L^\Phi(L) = L^\phi(L).
\]
It is clear that $\aff\dom \tilde h= L$.
\end{proof}

Taking a small enough space and a strong enough topology, it may be possible to establish relative continuity of $Eh$ throughout $\rcore\dom Eh$ even when the integrand involves nonaffine pointwise constraints. The set-valued mappings
\[
\omega\mapsto\dom h(\cdot,\omega),\quad \omega\mapsto\aff\dom h(\cdot,\omega)
\]
are measurable by \cite[Proposition 14.28 and Exercise~14.12]{rw98}.

The following is obtained by modifying arguments in the proof of \cite[Theorem~2]{roc71} which required, in particular, that $\aff\dom h=\reals^n$ almost surely; see also \cite[Theorem~4]{pp18a}.

\begin{example}\thlabel{thm:rinte}
Assume that $Eh:L^\infty\to\ereals$ is lsc and finite on
\[
\D := \{u\in L^\infty(\dom h) \mid \exists r>0:\ \uball_r(u)\cap\aff \dom h\subseteq\dom h\ \text{a.s.}\}.
\]
If $\D\ne\emptyset$, then $\D\subset\rcore\dom Eh$,
\[
\aff\dom Eh=L^\infty(\aff\dom h)
\]
and $Eh$ is strongly relatively continuous throughout $\rcore\dom Eh$.
\end{example}

\begin{proof}
Let $u\in\D$. There is an $r>0$ such that $\uball_r(u)\cap\aff \dom h\subseteq\dom h$ almost surely. Given $u'\in L^\infty(\aff\dom h)$, there is a $\lambda>0$ such that 
\[
\uball_{r/2}(\lambda(u'-u)+u)\subset\uball_r(u).
\]
Since $\lambda(u'-u)+u\in\aff\dom h$, we get $\lambda (u'-u)+u\in\D$ and thus, $L^\infty(\aff\dom h-u)\subseteq\pos(\D-u)$. By assumption,
\[
\D\subseteq\dom Eh\subseteq L^\infty(\dom h)\subseteq L^\infty(\aff\dom h),
\]
so
\[
\pos(\D-u) = \pos(\dom Eh-u) =\pos(L^\infty(\dom h)-u) = L^\infty(\aff\dom h)-u.
\]
Since $L^\infty(\aff\dom h)$ is affine, the first two claims follow from \thref{lem:rcore}. Since $L^\infty(\aff\dom h)$ is closed, the last claim follows from \thref{thm:Ehcont}.
\end{proof}


The following gives a simple example where $\pos\dom Eh$ is linear but not closed. In other words, we get $0\in\rcore\dom Eh$ but $\aff\dom Eh$ is not closed. \thref{cor:relcontiff} thus implies that the function $Eh$ is not relatively continuous on $\rcore\dom Eh$.

\begin{cexample}
Let $\U=L^\infty$, $h(u,\omega)=\delta_S(u,\omega)$, where $S(\omega):=\{u\mid |u|\le \eta(\omega)\}$ for strictly positive $\eta\in L^\infty$ such that $1/\eta\notin L^\infty$. Then
\[
\pos \dom Eh=\{u\in L^\infty\mid \exists \lambda>0: |u|\le \lambda\eta\}
\]
is linear but not closed in $L^\infty$. Indeed, defining $u^\nu:=\min\{\sqrt{\eta},\nu\eta\}$, we have $u^\nu\in\nu\dom Eh\subset\pos\dom Eh$ for all $\nu$ and $u^\nu\to\sqrt{\eta}$ in $L^\infty$ but, since $1/\eta\notin L^\infty$, we have $\sqrt{\eta}\notin\nu\dom Ef$ for all $\nu$.
\end{cexample}


In problems with nonlinear pointwise constraints, the strong relative core tends to be empty except in the strong topology of $L^\infty$.

\begin{example}\thlabel{ex:UneLinfty}
Let $S$ be a closed convex random set and assume that the probability space is atomless and that, for every $w\in L^\infty$ and $(A^\nu)_{\nu=1}^\infty\subset \F$ with $P(A^\nu)\searrow 0$, $1_{A^\nu} w\to 0$ in the strong topology of $\U$. Then 
$u\in\rcore\U(S)$ with $\aff\U(S)$ closed if and only if $u\in S$ almost surely and $S$ is affine-valued.
\end{example}

\begin{proof}
Sufficiency is clear. To prove necessity, let $u\in\rcore\U(S)$, $\aff\U(S)$ be closed and $S_u:=S-u$. By \thref{fb}, the set $\U(S_u)$ is a neighborhood of the origin in $\aff\U(S_u)$. Assume, for a contradiction, that $S_u$ is not linear-valued. 
Then $\bar A:=\{\omega\mid S_u(\cdot,\omega) \neq \aff S_u(\omega)\}$ is not a null set. Let $w\in L^0$ be such that $w\in S_u\setminus\rinte S_u$ on $\bar A$. There exists constants $r<R$ such that $A:=\{\omega\in \bar A\mid r<|w(\omega)|<R\}$ is not a null set. Now $1_Aw \in \U(S)$, so $1_A 2 w \in \aff\U(S)$ but $1_A 2w\notin S_u$ almost surely. Since $\U(S-u)$ is a neighborhood of the origin in $\aff\U(S-u)$, there exists, by the assumption in the statement, a non-null $A'\subseteq A$ such that  $2w1_{A'}\in\U(S-u)$. Thus, $2w1_{A'}\in S_u$ almost surely, which is a contradiction. 
\end{proof}



The topological assumption in \thref{ex:UneLinfty} holds, in particular, if any sequence $(u^\nu)_{\nu=1}^\infty$ in $L^\infty$ with $|u^\nu|\downto 0$ converges strongly to zero in $\U$. This condition fails in the strong topology of $L^\infty$ but holds e.g.\ in $L^p$ spaces with $p>1$ or, more generally, in Orlicz spaces associated with finite Young functions $\Phi$. These and more examples can be found in \cite{pp220}.

\subsection{The strong dual}\label{sec:sdU}

As in \cite{pp22b}, we assume that $\U$ and $\Y$ are solid decomposable spaces of random variables in separating duality under the bilinear form $E[u\cdot y]$. We will also assume that $\U$ is Fr\'echet under a given strong topology and that the corresponding topological dual $\U^*$ can be expressed as
\[
\U^* = \Y\oplus \Y^s
\]
in the sense that, for every $y\in\U^*$, there exist unique $y^c\in\Y$ and $y^s\in\Y^s$ such that
\[
\langle u,y\rangle = E[u\cdot y^c] + \langle u,y^s\rangle
\]
for all $u\in\U$. Here and in what follows, $\Y^s$ denotes the elements $y\in\U^*$ that are {\em singular} in the sense that, for every $u\in\U$, there exists a decreasing sequence $(A^\nu)_{\nu=1}^\infty\subset\F$ such that $P(A^\nu)\searrow 0$ and 
\[
\langle 1_{\Omega\backslash A^\nu} u,y\rangle= 0\quad\forall \nu=1,2,\ldots.
\]
By \thref{lem:kothe}, $\Y$ is necessarily the K\"othe dual of $\U$. Most familiar Banach spaces of random variables have a topological dual of the above form; see \cite{pp220}.

\begin{lemma}\thlabel{lem:EGcont}
  Given a $\sigma$-algebra  $\G\subset\F$ with $E^\G\U\subset\U$, the mapping $E^\G:\U\to\U$ is both strongly and $(\sigma(\U,\Y),\sigma(\U,\Y))$-continuous. Moreover, the adjoint $(E^\G)^*:\U^*\to\U^*$ is given by $(E^\G)^*y=E^\G y$ for every $y\in\Y$. In particular, $E^\G\Y\subset\Y$.
\end{lemma}

\begin{proof}
By \thref{lem:kothe}, $\Y$ is the K\"othe dual of $\U$. Thus, by \cite[Lemma~3]{pp22b}, $E^\G$ is continuous with respect to $\sigma(\U,\Y)$-topology. In particular, $E^\G$ has $\sigma(\U,\Y)\times\sigma(\U,\Y)$-closed graph. Thus the graph is strongly closed as well. Since $\U$ is Fr\'echet, the strong continuity follows from  the closed graph theorem. Given $y\in\Y$, \cite[Lemma~3]{pp22b} gives $E^\G y\in\Y$ and
\[
\langle u,(E^\G)^*y\rangle = \langle E^\G u,y\rangle = \langle u,E^\G y\rangle \quad
\]
for every $u\in\U$. Since $\U$ separates points in $\U^*$, we have  $(E^\G)^* y=E^\G y$.
\end{proof}

From now on, we will assume that $\X$ has the same properties as $\U$. That is, we assume that there is a topology under which $\X$ is a Fr\'echet space whose topological dual $\X^*$ can be expressed as
\[
\X^* = \V\oplus \V^s
\]
in the sense that, for every $v\in\X^*$, there exist unique $v^c\in\V$ and $v^s\in\V^s$ such that, for all $x\in\X$
\[
\langle x,v\rangle = E[x\cdot v^c] + \langle x,v^s\rangle.
\]
Here $\V^s$ is the set of singular elements of $\X^*$.

\begin{lemma}\thlabel{lem:A}
Let $A$ be a random matrix such that $A x\in\U$ for all $x\in\X$. Then $A^*y\in\V$ for all $y\in\Y$ and the linear mapping $\A:\X\to\U$ defined pointwise by
\[
\A x=Ax\quad a.s.
\] 
is both strongly and $(\sigma(\X,\V),\sigma(\U,\Y))$-continuous. Moreover, the adjoint $\A^*:\U^*\to\X^*$ is given by $\A^* y= A^*y$ for every $y\in\Y$.
\end{lemma}

\begin{proof}
By \thref{lem:kothe}, $\V$ is the K\"othe dual of $\X$ so the $(\sigma(\X,\V),\sigma(\U,\Y))$-continuity follows from \cite[Lemma~4]{pp22b}. In particular, $\gph\A$ is weakly closed. Weak closedness implies strong closedness. Since $\U$ and $\X$ are Fr\'echet, the strong continuity follows from  the closed graph theorem. Given $y\in\Y$, \cite[Lemma~4]{pp22b} gives $A^*y\in\V$ and
\[
\langle x,\A^*y\rangle = \langle Au,y\rangle = \langle u, A^*y\rangle
\]
for every $u\in\U$. Since $\U$ separates points in $\U^*$, we have  $\A^* y=A^*y$.
\end{proof}

The following two examples are from \cite{pp220}.

\begin{example}[Lebesgue spaces]\thlabel{ex:Lp5}
Let $\U=L^p$ and $\Y=L^q$ be the usual Lebesgue spaces, where $p$ and $q$ are conjugate exponents, and let the strong topology be the $L^p$-norm topology. Then $\U$ is Banach and its dual can be expressed as
\[
\U^*=\Y\oplus \Y^s.
\]
If $p<\infty$, then $\Y^s=\{0\}$. If $p=\infty$ and $(\Omega,\F,P)$ is atomless, then $\Y^s\ne\{0\}$.
\end{example}

\begin{example}[Musielak-Orlicz spaces]\thlabel{ex:mo5}
Let $\U=L^\Phi$ and $\Y=L^{\Phi^*}$ be as in \thref{ex:orliczcont} and let the strong topology be the topology generated by the norm
\[
\|u\|_{L^\Phi}:=\inf_{\alpha\in\reals_+}\{\alpha\mid E\Phi(u/\alpha)\le 1\}
\]
Then $\U$ is Banach and its dual can be expressed as
\[
\U^*=\Y\oplus \Y^s.
\]
If $\dom E\Phi$ is a cone, then $\Y^s=\{0\}$. We have $\Y^s\ne\{0\}$ if $L^\infty$ is not dense in $L^\Phi$ or if $P$ is atomless and $\Phi(a,\cdot)\notin L^1$ for some $a>0$.

If $\Phi(a,\cdot)\in L^1$ for all $a>0$, then the Morse heart 
\[
M^\Phi :=\{ u\in L^0 \mid E\Phi(\alpha u)<\infty\ \forall \alpha\in\reals_+ \}
\]
is Banach under the relative topology of $L^\Phi$ and $L^{\Phi^*}$ is its strong dual.
\end{example}

Much like Musielak-Orlicz spaces are often chosen so that a given objective is continuous (see \thref{ex:orliczcont}), the following example shows how to design the space $\U$ so that pointwise linear constraints define a continuous surjection to $\U$. Given a random matrix $A\in L^0(\reals^{m\times n})$ with $\rge A=\reals^m$ almost surely, its scenario-wise Moore-Penrose inverse 
\[
A^\dagger:=A^*(A^*A)^{-1}
\]
is measurable.

\begin{example}[Range spaces] \thlabel{ex:scaled}
 Let $A\in L^0(\reals^{m\times n})$,
\begin{align*}
\U &:=\{Ax\mid x\in\X\},\\
\Y &:=\{y\in L^0 \mid A^*y\in\V\}
\end{align*}
and assume that $L^\infty\subset\U$ and $L^\infty\subset\Y$. We endow $\U$ with the ``final topology'' induced by $A$, i.e.\ the strongest locally convex topology under which the linear mapping $x\mapsto Ax$ is continuous. Then $\U$ is a decomposable Fr\'echet space.  If $\X$ is such that
\[
\X=\{x\in L^0\mid |x|\in\X_0\},
\]
where $\X_0$ is a decomposable solid space of scalar random variables, then
\[
\U^*=\Y\oplus\Y^s,
\]
where $\Y^s$ is the set of singular elements of $\U^*$. If, in addition, $|A||A^\dagger|\in L^\infty$ then $\U$ is solid and
\[
\U=\{A_1 x\mid x\in\X\}\times\dots\times\{A_m x\mid x\in\X\},
\]
where $A_i$ is the $i$-th row of $A$.
\end{example}

\begin{proof}
We apply \thref{finalU} to the linear mapping $\A:\X\to L^0(\reals^m)$ defined pointwise by $(\A x)(\omega):=A(\omega)x(\omega)$. Clearly $\U=\rge\A$. We have $\ker\A=\{x\in\X\mid Ax=0\}$ which is weakly closed, by \cite[Corollary~7]{pp22b}, and thus, strongly closed. By \thref{finalU}, $\U$ is then Fr\'echet and, for every $u^*\in\U^*$, there is a unique $v\in(\ker\A)^\perp$ such that
\[
\langle \A x,u^*\rangle = \langle x,v\rangle \quad\forall x\in\X.
\]

Every $u\in\rge\A$ can be expressed as $u=\A(A^\dagger u)$, where $A^\dagger x$ is an element of $\X$. Indeed, $u=\A x$ for some $x\in\X$ so
\[
A^\dagger u = A^*(A^*A)^{-1}Ax
\]
which is the scenario-wise Euclidean projection of $x(\omega)$ on $\rge A^*$. Thus, $|A^\dagger Ax|\le|x|$ almost surely. Since $\X=\{x\in L^0\mid |x|\in\X_0\}$ , we get $A^\dagger u\in\X$ for all $u\in\rge\A$. It follows that, for every $u^*\in(\rge\A)^*$, there is a $v\in(\ker\A)^\perp$ such that
\begin{align*}
  \langle u,u^*\rangle &= \langle A^\dagger u,v\rangle\\
  &=E[A^\dagger u\cdot v^c]+\langle A^\dagger u,v^s\rangle\\
  &=E[ u\cdot (A^\dagger)^*v^c]+\langle A^\dagger u,v^s\rangle
\end{align*}
for all $u\in\rge\A$. Here $(A^\dagger)^* v^c\in\Y$, since $A^*(A^\dagger)^* v^c =v^c \in \V$. Moreover, $u\mapsto \langle A^\dagger u,v^s\rangle$ is singular, which follows directly from the definition. Thus every continuous linear functional can be expressed by an element of $\Y\oplus\Y^s$.

Conversely, let $(y^c,y^s)\in\Y\oplus\Y^s$. By definition of the final topology, $u\mapsto E[u\cdot y^c]$ is continuous on $\U$ if and only if $x\mapsto E[(Ax)\cdot y^c]$ is continuous on $\X$. This holds by definition of $\Y$. The singular component $y^s$ is continuous by definition.

Assume now that $|A||A^\dagger|\in L^\infty$. Let $u\in\U$ and  $u' \in L^0(\reals^m)$ with  $|u'|\le |u|$. There exists $x\in\X$ with $Ax=u$, so 
\[
|A^\dagger u'|\le |A^\dagger| |u| \le |A^\dagger||A||x|.
\]
The right side belongs to $\X_0$, so  $A^\dagger u'\in\X$ by assumption. Thus $A^\dagger u'\in\X$, so $u'=A A^\dagger u' \in\U$ which shows that $\U$ is solid. The last claim follows from the fact that solid spaces are Cartesian products.
\end{proof}

\subsection{Conjugates of integral functionals}

Let $\U$ and $\U^*$ be as in Section~\ref{sec:sdU}. This section computes conjugates and subdifferentials of integral functionals with respect to the duality pairing of $\U$ and $\U^*$. The main result of this section makes use of the following simple observation.

\begin{lemma}\thlabel{usac}
We have 
\[
\limsup_{\nu\to\infty} Eh(1_{A^\nu}u+1_{\Omega\backslash A^\nu}\bar u )\le Eh(\bar u)
\]
for all $u,\bar u\in\dom Eh$ and decreasing $(A^\nu)_{\nu=1}^\infty\subset\F$ with $P(A^\nu)\downto 0$.
\end{lemma}
\begin{proof}
For any $\bar u,u\in \dom Eh$ and $A\in\F$, the right side in
\[
h(1_{A} u+ 1_{\Omega\backslash A} \bar u)\le \max \{ h(u),h(\bar u)\}
\]
is integrable, so the claim follows from Fatou's lemma.
\end{proof}

The following extends \cite[Theorem~10]{pp18a} which in turn refines \cite[Corollary~1B]{roc71}. Both \cite{roc71} and \cite{pp18a} studied integral functionals on the space $L^\infty$ of essentially bounded random variables. The following extends \cite[Theorem~10]{pp18a} to general Fr\'echet spaces that satisfy the assumptions Section~\ref{sec:sdU}. It also yields \cite[Theorem~2.6]{koz80} in the case of Orlicz spaces of $\reals^n$-valued random variables. 

\begin{lemma}\thlabel{lem:usac}
Let $Eh$ be proper. If $\bar u\in\U$, $y\in\U^*$ and $\epsilon\ge 0$ are such that
\begin{equation}\label{esg}
Eh(u)\ge Eh(\bar u) + \langle u-\bar u,y\rangle-\epsilon\quad \forall u\in \U,
\end{equation}
then
\begin{equation}\label{esgac}
Eh(u)\ge Eh(\bar u) + \langle u-\bar u,y^c\rangle-\epsilon\quad \forall u\in\U,
\end{equation}
and
\begin{equation}\label{esgas}
0\ge\langle u-\bar u,y^s\rangle-\epsilon\quad\forall u\in\dom Eh.
\end{equation}
\end{lemma}

\begin{proof}
Let $u\in\dom Eh$ and let $(A^\nu)_{\nu=1}^\infty\subset\F$ be the decreasing sequence of sets in the characterization of the singular component $y^s$ of $y$. Let $u^\nu:= \one_{A^\nu}\bar u+\one_{\Omega\setminus A^\nu}u$. \thref{usac} and \eqref{esg} give
\[
Eh(u)\ge\limsup Eh(u^\nu)\ge Eh(\bar u)+\limsup\langle u^\nu-\bar u, y\rangle - \epsilon.
\]
Since $u^\nu-\bar u= \one_{\Omega\setminus A^\nu}(u-\bar u)$, we have
\[
\langle u^\nu-\bar u,y\rangle = \langle u^\nu-\bar u,y^c\rangle \to \langle u-\bar u,y^c\rangle
\]
so \eqref{esgac} holds.

Now let $u^\nu:= \one_{A^\nu}{u} + \one_{\Omega\setminus A^\nu}\bar u$ so that
\[
Eh(\bar u)\ge\limsup Eh(u^\nu)\ge Eh(\bar u)+\limsup\langle u^\nu-\bar u,y\rangle - \epsilon.
\]
Since $u^\nu-\bar u=\one_{A^\nu}(u-\bar u)$, we have
\begin{align*}
\langle u^\nu-\bar u,y\rangle = \langle u^\nu-\bar u,y^c\rangle+\langle u^\nu-\bar u,y^s\rangle\to \langle u-\bar u,y^s\rangle
\end{align*}
so $\langle u-\bar u,y^s\rangle - \epsilon\le 0$. Since $u\in\dom Eh$ was arbitrary, \eqref{esgas} holds.
\end{proof}

Recall that $\U^*=\Y\oplus\Y^s$. We now come to the main result of this section which states that, at each $y=y^c+y^s\in\U^*$, the conjugate of $Eh$ equals the sum of the integral functional of $h^*$ evaluated at $y^c$ and the recession function of $(Eh)^*$ evaluated at $y^s$; see \thref{domrec}. The following extends the main results of \cite{roc71} on integral functionals on $L^\infty$ to more general Fr\'echet spaces that satisfy the assumptions of Section~\ref{sec:sdU}.

\begin{theorem}\thlabel{thm:usac}
If $Eh$ is proper, then
\begin{enumerate}
\item
  $(Eh)^*(y) = Eh^*(y^c)+\sigma_{\dom Eh}(y^s)\quad\forall y\in\U^*$,
\item
  $Eh$ is strongly lsc if and only if it is $\sigma(\U,\Y)$-lsc,
\item
  $y\in\partial Eh(\bar u)$ if and only if $y^c\in\partial Eh(\bar u)$ and $y^s\in N_{\dom Eh}(\bar u)$.
\end{enumerate}
\end{theorem}

\begin{proof}
We have
\begin{align*}
(Eh)^*(y) &=\sup_{u}\{\langle u,y\rangle - Eh(u)\}\\
&=\sup_{u}\{\langle u,y^c\rangle - Eh(u)+\langle u,y^s\rangle-\delta_{\dom Eh}(u)\}\\
&\le Eh^*(y^c)+\sigma_{\dom Eh}(y^s)
\end{align*}
while for every $y\in\dom (Eh)^*$ and $\epsilon>0$, there exists $\bar u\in\dom Eh$ such that \eqref{esg} holds. By \thref{lem:usac} and Fenchel's inequality, 
\begin{align*}
 (Eh)^*(y^c)+\sigma_{\dom Eh}(y^s) \le \langle \bar u,y^c\rangle -Eh(\bar u)+\epsilon + \langle \bar u, y^s\rangle +\epsilon \le (Eh)^*(y)+2\epsilon.
\end{align*}
This gives the expression for the conjugate since $\epsilon>0$ was arbitrary.

It is clear that if $Eh$ is $\sigma(\U,\Y)$-lsc it is strongly lsc. When $Eh$ is proper and strongly lsc, the biconjugate theorem says that it has a proper conjugate. The first claim then implies that $Eh^*$ is proper on $\Y$. The second claim thus follows from \cite[Corollary~6]{pp220}.

When $\epsilon=0$, \thref{lem:usac} says that if $y\in\partial Eh(\bar u)$ then $y^c\in\partial Eh(\bar u)$ and $y^s\in N_{\dom Eh}(\bar u)$. The converse implication follows simply by adding \eqref{esgac} and \eqref{esgas} together.
\end{proof}

\thref{thm:usac} implies, in particular, that strong subdifferentiability of $Eh$ implies subdifferentiability with respect to the pairing of $\U$ with $\Y$. 

\begin{corollary}\thlabel{cor:IFmc}
Assume that $Eh$ is closed proper and finite on $\U$. Then $Eh$ is continuous both strongly and in $\tau(\U,\Y)$.
\end{corollary}

\begin{proof}
Strong continuity follows from \thref{lscbar}. By Bourbaki-Alaoglu, level sets of $(Eh)^*$ are $\sigma(\U^*,\U)$-compact. By \thref{thm:usac}, the finiteness of $Eh$ implies that all the level sets of $(Eh)^*$ are contained in $\Y$, so they are $\sigma(\Y,\U)$-compact. By the converse of Bourbaki-Alaoglu, $Eh$ is $\tau(\U,\Y)$-continuous.
\end{proof}



\section{The strong dual problem}\label{sec:sdp}

From now on, we will assume that both $\X$ and $\U$ satisfy the assumptions of Section~\ref{sec:sdU}. More precisely, $\X$ and $\U$ are endowed with a {\em strong} topology under which they are Fr\'echet, the dual of $\U^*$ can be identified with $\Y\oplus\Y^s$ and, similarly, the topological dual of $\X^*$ can be identified with $\V\oplus\V^s$ in the sense that, for each $v\in\X^*$, there exist unique $v^c\in\V$ and $v^s\in\V^s$ such that
\[
\langle x,v\rangle = E[x\cdot v^c] + \langle x,v^s\rangle
\]
for all $x\in\X$. Here, $\V^s$ denotes the singular elements of $\X^*$.

Much like in \cite{pp22b}, we apply the general conjugate duality framework of \cite{roc74} to \eqref{spx} with
\[
F(x,z,u) := Ef(x,u) + \delta_\N(x-z),
\]
but this time, with respect to the pairings of $\X$ and $\U$ with $\X^*$ and $\U^*$, respectively. This gives the {\em strong dual problem}
\begin{equation}\label{sd}\tag{$D_s$}
\maximize\quad \langle\bar  u,y\rangle - \varphi^*(p,y)\quad\ovr (p,y)\in\X^*\times\U^*,
\end{equation}
where $\varphi^*$ is the conjugate of $\varphi$ with respect the strong pairings. As long as the optimum value of \eqref{spx} is finite, the subgradients of the optimum value function $\varphi$ in the space $\X^*\times\U^*$ are solutions of \eqref{sd}; see \thref{sdvarphi}. Establishing the existence of a subgradient in the strong dual may be significantly easier than to establish the existence in the smaller space $\V\times\Y\subset\X^*\times\U^*$. Indeed, since $\X$ and $\U$ are Fr\'echet, \thref{thm:core} gives general conditions for the subdifferentiability. The special structure of the conjugates given in Section~\ref{sec:cif} will then allow us to show in Section~\ref{sec:recourse} that the existence of solutions to \eqref{sd} implies the existence of solutions of the original dual problem \eqref{d} and, moreover, that the optimum values are equal.

We will denote the orthogonal complement of $\X_a$ in the strong dual of $\X$ by
\[
\X_a^\circ:=\{p\in\X^*\mid \langle x,p\rangle=0 \ \forall x\in\X_a\}.
\]
Clearly $\X_a^\perp\subset\X_a^\circ$. The following lemma gives expressions for the conjugates of $F$ and the optimum value function
\[
\varphi(z,u):=\inf_{x\in\X}\{Ef(x,u)\,|\,x-z\in\N\}
\]
with respect to the strong pairings.

\begin{lemma}\thlabel{lem:sdual}
If $\dom Ef\cap(\X\times\U)\ne\emptyset$, then
\[
F^*(v,p,y) = Ef^*((v+p)^c,y^c)+\sigma_{\dom Ef}((v+p)^s,y^s)+\delta_{\X_a^\circ}(p)
\]
and, in particular, 
\[
\varphi^*(p,y) = Ef^*(p^c,y^c)+\sigma_{\dom Ef}(p^s,y^s)+\delta_{\X_a^\circ}(p).
\]
If, in addition, $\dom Ef^*\cap(\V\times\Y)\ne\emptyset$, then $F$ is proper and strongly lsc.
\end{lemma}

\begin{proof}
By \thref{thm:usac},
\begin{align*}
F^*(v,p,y) &= \sup_{x\in\X,z\in\X,u\in\U}\{\langle x,v\rangle + \langle z,p\rangle + \langle u,y\rangle - Ef(x,u) \,|\, x-z\in\X_a\}\\
&= \sup_{x\in \X,z'\in\X,u\in\U}\{\langle x,v+p\rangle + \langle u,y\rangle - Ef(x,u) - \langle z',p\rangle] \,|\, z'\in\X_a\}\\
&=Ef^*((v+p)^c,y^c)+\sigma_{\dom Ef}((v+p)^s,y^s)+\delta_{\X_a^\circ}(p).
\end{align*}
When $\dom Ef^*\ne\emptyset$, \cite[Lemma~13]{pp22b} says that $Ef$ is $\sigma(\X,\V)\times\sigma(\U,\Y)$-lsc and proper. The last claim thus follows from the fact that $F$ is the sum of lsc functions. Clearly, weak lower semicontinuity implies strong semicontinuity and $F$ is proper when $Ef$ is proper.
\end{proof}

As an immediate corollary, we get the following.

\begin{theorem}\thlabel{thm:dualproblem}
If $\dom Ef\cap(\X\times\U)\ne\emptyset$, the strong dual problem \eqref{sd} can be written as
\begin{equation}\label{sd}\tag{$D_s$}
\maximize\ \langle \bar u,y\rangle -  Ef^*(p^c,y^c)-\sigma_{\dom Ef}(p^s,y^s)\ \ovr\ (p,y)\in\X_a^\circ\times\U^*.
\end{equation}
\end{theorem}

Restricting the dual variables $(p,y)$ to $\V\times\Y$, problem \eqref{sd} reduces to \eqref{d}, so
\[
\inf\eqref{sd}\le \inf\eqref{d}.
\]
In particular, if there is no duality gap between \eqref{spx} and \eqref{d}, the same holds between \eqref{spx} and \eqref{sd}. Clearly, the two dual problems coincide if $\Y^s=\{0\}$ and $\X^s=\{0\}$ as happens e.g.\ when $\U=\X=L^p$ with $p<\infty$; see \thref{ex:Lp5}. The problems coincide also in the special case where $\dom Ef=\X\times\U$ since then, $\sigma_{\dom Ef}=\delta_{\{(0,0)\}}$. Section~\ref{sec:recourse} below gives more general conditions under which the optimum values of \eqref{d} and \eqref{sd} coincide and one has a solution if and only if the other one does.

The following gives sufficient conditions for subdifferentiability of $\varphi$ with respect to $\X^*\times\U^*$ and thus, by \cite[Theorem~16]{roc74}, for the absence of a duality gap between \eqref{spx} and \eqref{sd} and for the existence of solutions for the latter.

\begin{theorem}\thlabel{thm:ssd}
If $Ef$ is proper and strongly closed and
\[
\pos(\dom\varphi-(0,\bar u))
\]
is linear and closed, then the strong dual problem \eqref{sd} has a solution and
\[
\inf\eqref{spx}=\sup\eqref{sd}.
\]\end{theorem}

\begin{proof}
Follows directly from \thref{thm:core} and \cite[Theorem~16]{roc74}.
\end{proof}

The set $\pos(\dom\varphi-(0,\bar u))$ is linear and strongly closed if $(0,\bar u)$ belongs to the strong interior of $\dom\varphi$. A simple extension of this is to require that the affine hull of $\dom\varphi$ is strongly closed and $(0,\bar u)$ belongs to the corresponding relative interior of $\dom\varphi$. The following lemma will be useful in providing more general conditions.

\begin{lemma}\thlabel{lem:posdomphi}
If $\bar x\in\X_a$ is feasible, then
\[
\pos(\dom\varphi-(0,\bar u)) = \pos(\dom Ef-(\bar x,\bar u)) - \X_a\times\{0\}
\]
\end{lemma}

\begin{proof}
By definition,
\[
\dom\varphi = \dom Ef - \X_a\times\{0\}
\]
so
\[
\dom\varphi-(0,\bar u) = \dom Ef -(\bar x,\bar u)- \X_a\times\{0\}
\]
for any $\bar x\in\X_a$. If $\bar x\in\X_a$ is feasible, then $(\bar x,\bar u)\in\dom Ef$ and
\[
\pos(\dom\varphi-(0,\bar u)) = \pos(\dom Ef -(\bar x,\bar u)) - \X_a\times\{0\},
\]
by \thref{lem:hulls}.
\end{proof}

We say that an $\bar x\in\X_a$ is {\em strictly feasible} if
\[
\pos(\dom Ef-(\bar x,\bar u))
\]
is linear and closed. If such a point exists, we say that the problem \eqref{spx} is {\em strictly feasible}. This happens, in particular, if there exists an $\bar x\in\X_a$ such that $(\bar x,\bar u)\in\inte\dom Ef$. This corresponds to the classical {\em Slater condition} in optimization theory. The condition implies that $\pos(\dom Ef-(\bar x,\bar u))=\X\times\U$ and thus, by \thref{lem:posdomphi}, $\pos(\dom\varphi-\bar u)=\X\times\U$, so $\varphi$ is strongly subdifferentiable, by \thref{thm:ssd}. More generally, if there exists a strictly feasible $\bar x$ such that $\X_a\times\{0\}\subset\pos(\dom Ef-(\bar x,\bar u))$, we have $\pos(\dom\varphi-\bar u)=\pos(\dom Ef-(\bar x,\bar u))$ and, by \thref{thm:ssd} again, $\varphi$ is subdifferentiable.

In general, \thref{lem:posdomphi} reduces the condition in \thref{thm:ssd} to the question on the closedness of the sum of two sets. Under strict feasibility, both sets are linear and closed. \thref{thm:rinte} gives the following sufficient condition for strict feasibility when $\X=L^\infty$ and $\U=L^\infty$.

\begin{example}
Let $\X=L^\infty$ and $\U=L^\infty$ and assume that $Ef$ is finite on the set
\[
\D := \{(x,u)\in(\X\times\U)(\dom f) \mid \exists r>0:\ \uball_r(x,u)\cap\aff \dom f\subseteq\dom f\ P\text{-a.e.}\}.
\]
If $\bar x\in\X_a$ is such that $(\bar x,\bar u)\in\D$, then $\bar x$ is strictly feasible and
\[
\pos(\dom Ef-(\bar x,\bar u)) = \{(x,u)\in\X\times\U\mid (x,u)\in\aff\dom f-(\bar x,\bar u)\ a.s.\}.
\]
\end{example}

Closedness of sums of closed linear subspaces is a well-studied question in functional analysis for which various sufficient as well as necessary conditions have been given. In many applications, of stochastic optimization, the problem has additional structure that can be used to establish the closedness of the sum in \thref{lem:posdomphi}. The following lemma gives a sufficient condition that may look awkward at first glance but turns out to be convenient in many applications.

\begin{lemma}\thlabel{lem:pircore}
Let $\bar x\in\X_a$ be strictly feasible and let $\L:=\pos(\dom Ef-(\bar x,\bar u))$. If there exists a strongly continuous linear idempotent mapping $\pi$ on $\X\times\U$ such that $\pi\L\subseteq\L$ and
\[
\L-\X_a\times\{0\}=\L+\rge\pi,
\]
then $\pos(\dom\varphi-(0,\bar u))$ is linear and closed. 
The equality holds, in particular, if
\[
\X_a\times\{0\}\subseteq\rge\pi\subseteq\L-\X_a\times\{0\}.
\]
\end{lemma}

\begin{proof}
Strict feasibility means that there exists $\bar x\in\X_a$ such that $\pos(\dom Ef-(\bar x,\bar u))=\L$ and that $\L$ is closed. 
By \thref{lem:posdomphi},
\[
\pos(\dom\varphi-(0,\bar u)) = \L - \X_a\times\{0\},
\]
so the first claim follows from \thref{lem:sss}. The additional inclusions give
\[
\L-\X_a\times\{0\}\subseteq\L+\rge\pi\subseteq\L-\X_a\times\{0\},
\]
which proves the second claim.
\end{proof}

In some applications, it may be difficult to find a mapping $\pi$ that satisfies the assumptions of \thref{lem:pircore}. In such situations, the formulation in \thref{ex:rcoregA} below may be more convenient. 

\begin{example}\thlabel{ex:rcoregA}
Let
\[
f(x,u,\omega) = g(x,A(\omega)x+u,\omega),
\]
where $g$ is a convex normal integrand on $\reals^n\times\reals^m\times\Omega$ and $A$ is a random matrix with $A\X\subseteq\U$. Assume that $\bar x\in\X_a$ is such that
\[
\K:=\pos(\dom Eg-(\bar x,A\bar x+\bar u))
\]
is linear and closed. Then $\bar x$ is strictly feasible. If there exists a strongly continuous linear idempotent mapping $\pi'$ on $\X\times\U$ such that $\pi'\K\subseteq\K$ and
\[
\K-(I,A)\X_a=\K+\rge\pi',
\]
then $\pos(\dom\varphi-(0,\bar u))$ is linear and closed. The equality holds, in particular, if
\[
(I,A)\X_a\subseteq\rge\pi'\subseteq\K-(I,A)\X_a.
\]
\end{example}

\begin{proof}
We have $\dom Ef=\{(x,u)\in\X\times\U\mid (x,Ax+u)\in\dom Eg\}=\bar A^{-1}\dom g$, where 
\[
\bar A = \begin{bmatrix}I & 0\\A & I\end{bmatrix}.
\]
By \thref{lem:hulls},
\[
\pos(\dom Ef-(\bar x,\bar u))=\bar A^{-1}\pos(\dom Eg-\bar A(\bar x,\bar u)) = \bar A^{-1}\K.
\]
Since $\K$ is linear and closed, by assumption, $\bar x$ is strictly feasible.

To complete the proof, we apply \thref{lem:pircore} with the idempotent mapping $\pi=\bar A^{-1}\pi'\bar A$. We have $\L=\bar A^{-1}\K$ so $\pi'\K\subseteq\K$ implies $\pi\L\subseteq\L$. Moreover,
\begin{align*}
  \L-\X_a\times\{0\} &= \{(x,u)\mid (x,Ax+u)\in \K\}-\X_a\times\{0\}\\
  &= \{(x,u)\mid (x,Ax+u)\in\K-(I,A)\X_a\}\\
  &= \bar A^{-1}(\K-(I,A)\X_a)
\end{align*}
and $\rge\pi = \bar A^{-1}\rge\pi'$, so the equality $\K-(I,A)\X_a=\K-\rge\pi'$ implies the equality in \thref{lem:pircore}. The additional inclusions give
\[
\K-(I,A)\X_a\subseteq\K+\rge\pi'\subseteq\K-(I,A)\X_a,
\]
which proves the last claim.
\end{proof}

\section{Existence of dual solutions}\label{sec:recourse}

The set of solutions of \eqref{d} is the intersection of the solution set of \eqref{sd} with $\V\times\Y$. In general, it may happen that \eqref{d} does not admit solutions even if \eqref{sd} does and there is no duality gap; see \cite[Example~2]{rw76b}. This section gives conditions under which, for each point $(p,y)\in\X^*\times\U^*$, there is a dual feasible $(\tilde p,\tilde y)\in\V\times\Y$ that achieves a dual objective value at least as good as $(p,y)$. Combined with the existence results for strong dual solutions in the previous section, we then obtain existence results for \eqref{d}.

Let $\X$ and $\U$ be as in Section~\ref{sec:sdp}. We will assume, in addition, that $E_t\X\subseteq\X$ for all $t$. By \thref{lem:EGcont}, this implies that the conditional expectations $E_t$ are both strongly and weakly continuous and that, when restricted to $\V$, the adjoints $E_t^*:\X^*\to\X^*$ coincide with the conditional expectations $E_t$. By \cite[Lemma~3]{pp22b} and \thref{lem:kothe}, we also have $E_t\V\subseteq\V$. Most familiar spaces of random variables satisfy the above assumptions.

We will denote the projection of the set $\dom Ef\subset\X\times\U$ on $\U$ by
\[
\dom_u Ef :=\{u\in\U\mid\exists\ x\in\X:\ (x,u)\in\dom Ef\}.
\]

\begin{assumption}\thlabel{ass:rec}
\mbox{}
\begin{enumerate}
\item
  $\pos(\dom_u Ef-\bar u)$ is linear and closed,
\item
  for every $t$ and $z\in\X$ with $(z,\bar u)\in \dom Ef$, there exists $\bar z\in\X$ with $(\bar z,\bar u)\in\dom Ef$ and $\bar z^t=E_t z^t$.
  \end{enumerate}
\end{assumption}

\thref{ass:rec} clearly holds if $\dom Ef=\X\times\U$. Note that, in this case, the support function of $\dom Ef$ is the indicator of the origin so feasible solutions in \eqref{sd} have $(p^s,u^s)=(0,0)$. More interesting sufficient conditions for \thref{ass:rec} are given in \thref{lem:rec0,ex:recinfty} below.

First, however, we will show that under \thref{ass:rec}, one can restrict dual variables to $\V\times\Y$ without worsening the objective value in the strong dual problem \eqref{sd}. The argument is based on a recursive application of the following somewhat technical lemma, the proof of which is an application of conjugate duality; see the appendix.

\begin{lemma}\thlabel{lem:spi}
\thref{ass:rec} implies that, for every $t$ and $(p,y)\in\X^*\times\U^*$ with $p_{t'}=0$ for $t'>t$, there exists a $\hat y\in\U^*$ such that
\[
\sigma_{\dom Ef}(E_t^*p,\hat y)-\langle\bar u,\hat y\rangle \le \sigma_{\dom Ef}(p,y)-\langle\bar u,y\rangle. 
\]
\end{lemma}

\begin{proof}
Define a convex function $g_t$ on $\X$ by
\begin{align*}
g_t(z) &:=\inf_{\tilde z\in\L_t}Ef(z+\tilde z,\bar u),
\end{align*}
where $\L_t:=\{z\in \X\mid z^t=0\}$. \thref{ass:rec}.2 means that
\[
\delta_{\dom g_t}(E_tz)\le \delta_{\dom g_t}(z)\quad\forall z\in\X.
\]
By \thref{lem:idemp,domrec}, this implies
\begin{equation}\label{eq:spi}
(g_t^*)^\infty(E_t^*p)\le (g_t^*)^\infty(p)\quad\forall p\in\X^*,
\end{equation}
A direct calculation gives
\[
g_t^* = (Ef(\cdot,\bar u))^*+\delta_{\L_t^\perp}(\cdot)
\]
where $\L_t^\perp=\{p\in\X^*\mid p_{t'}=0\ \forall t'>t\}$. Since the recession function of the sum is the sum of the recession functions, we get by \thref{domrec} that
\begin{align*}
(g_t^*)^\infty &=\delta_{\dom Ef}(\cdot,\bar u)^*+\delta_{\L_t^\perp}(\cdot),
\end{align*}
By Theorem~\ref{domcd}, \thref{ass:rec}.1 implies 
\[
\delta_{\dom Ef}(\cdot,\bar u)^*(p) = \inf_{y\in\U^*}\{\sigma_{\dom Ef}(p,y)-\langle\bar u,y\rangle\}
\]
where the infimum is attained for every $p\in\X^*$.  Combining this with the expression for $(g_t^*)^\infty$ and with \eqref{eq:spi} proves the claim.
\end{proof}

The following is the main result of this section. Its proof is based on \thref{lem:sdual,lem:spi}, The argument extends that of \cite[Theorem~3]{pp18a} which in turn simplifies the arguments of \cite{rw76}.

\begin{theorem}\thlabel{thm:spiu}
Under \thref{ass:rec}, there exists, for every $(p,y)\in\X_a^\circ\times\U^*$, a dual feasible $(\tilde p,\tilde y)\in\X_a^\perp\times\Y$ with
\begin{equation}\label{spiineq}
\varphi^*(\tilde p,\tilde y)-\langle \bar u,\tilde y\rangle \le \varphi^*(p,y)-\langle \bar u,y\rangle.
\end{equation}
\end{theorem}

\begin{proof}
Let $(p,y)\in\X_a^\circ\times\U^*$. We will show that if $p_{t'}^s=0$ for $t'>t$, then there exists $(\tilde p,\tilde y)\in\X_a^\circ\times\U^*$ satisfying \eqref{spiineq} and $\tilde p_{t'}^s=0$ for $t'\ge t$. Since $p^s_{t'}=0$ for $t'>T$, we then get, by induction, a $(\tilde p,\tilde y)\in\X_a^\perp\times\U^*$ satisfying \eqref{spiineq}.

By \thref{lem:sdual},
\[
\varphi^*(p,y)-\langle\bar u,y\rangle = Ef^*(p^c,y^c)-\langle\bar u,y^c\rangle+\sigma_{\dom Ef}(p^s,y^s)-\langle\bar u,y^s\rangle.
\]
Assume that $p^s_{t'}=0$ for $t'>t$. By \thref{lem:spi}, there is a $\hat y\in\U^*$ such that
\begin{align*}
\varphi^*(p,y)-\langle\bar u,y\rangle &\ge Ef^*(p^c,y^c)-\langle\bar u,y^c\rangle+\sigma_{\dom Ef}(E_t^*p^s,\hat y)-\langle\bar u,\hat y\rangle.
\end{align*}
Applying Fenchel's inequality to $Ef$ and the indicator of $\dom Ef$, now gives
\begin{align*}
\varphi^*(p,y)-\langle\bar u,y\rangle &\ge \sup_{x\in\X, u\in\U}\{\langle x,p^c+E_t^*p^s\rangle + \langle u,y^c+\hat y\rangle - Ef(x,u)\}-\langle\bar u,y^c+\hat y\rangle\\
&= Ef^*(\tilde p^c,\tilde y^c)+\sigma_{\dom Ef}(\tilde p^s,\tilde y^s)-\langle\bar u,\tilde y\rangle,
\end{align*}
where $\tilde p = p^c+E_t^*p^s$, $\tilde y = y^c+\hat y$ and the last equality holds by \thref{thm:usac}. Since, $p\in\X_a^\circ$, we have $E_t^*p_t=0$ so $E_t^*p_t^s=-E_t^*p_t^c=-E_tp_t^c$ and thus, $\tilde p^s_{t'}=0$ for every $t'\ge t$ as desired. It is easily checked that we still have $\tilde p\in\X_a^\circ$ so, by Theorem~\ref{lem:sdual} again,
\[
\varphi^*(\tilde p,\tilde y)-\langle\bar u,\tilde y\rangle\le\varphi^*(p,y)-\langle\bar u,y\rangle.
\]
Repeating the above procedure, we will get a $(\tilde p,\tilde y)\in\X_a^\perp\times\U^*$ satisfying \eqref{spiineq} as claimed. Since
\[
\sigma_{\dom Ef}(0,y)-\langle\bar u,y\rangle\ge 0\quad\forall y\in\U^*,
\]
we can also drop the singular part of $\tilde y$ with out increasing the value of $\varphi^*(\tilde p,\tilde y)-\langle\bar u,\tilde y\rangle$.
\end{proof}

It is clear that the optimum value of the strong dual minorizes that of the weak dual. \thref{thm:spiu} can thus be reformulated as follows.

\begin{theorem}\thlabel{thm:spi}
Under \thref{ass:rec}, the optimum values of \eqref{sd} and \eqref{d}  coincide and one has a solution if the other does. 
\end{theorem}

Combined with the results of Section~\ref{sec:sdp}, the above gives a two-step strategy for establishing the existence of solutions to \eqref{d}. First, one verifies subdifferentiability of the optimum value function $\varphi$ with respect to the pairing of $\X\times\U$ with $\X^*\times\U^*$ and second, one checks whether \thref{ass:rec} is satisfied. For strong subdifferentiability, we can use the sufficient conditions given in Section~\ref{sec:sdp}. Combining this with the optimality conditions from \cite{pp22b} gives the following.

\begin{theorem}\thlabel{thm:dexistence}
Assume that $\pos(\dom\varphi-(0,\bar u))$ is linear and closed and that \thref{ass:rec} holds. Then
\[
\inf\eqref{sp}=\inf\eqref{sp}=\sup\eqref{d}
\]
and the dual optimum is attained. Moreover, a primal feasible $x$ solves \eqref{sp} if and only if there exists a dual feasible $(p,y)$ with $(p,y)\in\partial f(x,\bar u)$ almost surely. 
\end{theorem}

\begin{proof}
By \thref{thm:ssd,thm:spi}, the assumptions imply the first two claims. The last claim follows from \cite[Corollary~27]{pp22b}.
\end{proof}




We end this section with some sufficient conditions for \thref{ass:rec}. Under strict feasibility, \thref{ass:rec}.1 can be stated as follows.

\begin{lemma}\thlabel{lem:rec0}
Let $\bar x\in\X_a$ be strictly feasible. \thref{ass:rec}.1 holds if and only if the projection of $\pos(\dom Ef-(\bar x,\bar u))$ to $\U$ is strongly closed.
\end{lemma}

\begin{proof}
\thref{ass:rec}.1 holds if and only if $\pos(\dom_uEf-\bar u)$ is linear and strongly closed. Strict feasibility means that there exists an $\bar x\in\X_a$ such that $\pos(\dom Ef-(\bar x,\bar u))$ is linear and strongly closed. We have
\[
\dom_uEf-\bar u = \{u\in\U\mid\exists x\in\X:\ (x,u)\in\dom Ef-(\bar x,\bar u)\}
\]
so, by \thref{lem:hulls},
\[
\pos(\dom_uEf-\bar u) = \{u\in\U\mid\exists x\in\X:\ (x,u)\in\pos(\dom Ef-(\bar x,\bar u))\}.
\]
This proves the claim.
\end{proof}

\begin{example}\thlabel{ex:rec1}
Consider \thref{ex:rcoregA} and assume that $\bar x\in\X_a$ is such that $(\bar x,A\bar x+\bar u)\in\core_s\dom Eg$. If $\{u\in\U\mid (0,u)\in\K-\gph A)$ is strongly closed, then \thref{ass:rec}.1 holds.
\end{example}

\begin{proof}
We have
\[
\dom Ef=\{(x,u)\in\X\times\U\mid (x,Ax+u)\in\dom Eg\}
\]
so
\[
\pos(\dom Ef-(\bar x,\bar u)) = \{(x,u)\in\X\times\U\mid (x,Ax+u)\in\pos(Eg-(\bar x,A\bar x,\bar u))\}.
\]
The projection of this set to $\U$ can be written as
\begin{align*}
  &\{u\in\U\mid\exists x\in\X: (x,u)\in\pos(\dom Ef-(\bar x,\bar u)\} \\
  &=\{u\in\U\mid \exists x\in\X:\ (x,Ax+u)\in\K\}\\
  &=\{u\in\U\mid (0,u)\in\K-\gph A\}
\end{align*}
so the claim follows from \thref{lem:rec0}.
\end{proof}



We say that a set-valued mapping $S:\Omega\tos\reals^{n_0}\times\cdots\times\reals^{n_T}$ is {\em $\FF$-adapted} if, for each $t=0,\ldots,T$, the projection
\[
S^t(\omega):=\{x^t \in\reals^{n^t}\,|\, x\in S(\omega)\}
\]
is $\F_t$-measurable. If $S(\omega)=\{s(\omega)\}$ for a stochastic process $s$, then $S$ is adapted if and only if the process $s$ is adapted.

Our condition in the following example extends the "bounded recourse condition" of \cite{rw83} from problems of Bolza to the general setting. If $S$ is essentially bounded, then the condition becomes the ``relatively complete recourse'' condition from \cite{rw76}.

\begin{example}[bounded recourse condition]\thlabel{ex:recinfty}
If $\X=L^\infty$ then \thref{ass:rec}.2 holds if, for all $r$ large enough, the mapping
\[
S_r(\omega):=\{x\in\uball_r\mid (x,\bar u(\omega))\in\dom f(\omega)\}
\]
is adapted and if $(x,\bar u)\in\dom Ef$ for all $x\in L^0(S_r)$.

If the set-valued mapping $S(\omega):=\{x\in\reals^n\mid (x,\bar u(\omega))\in\dom f(\omega)\}$ is essentially bounded, the above clearly holds if $S$ is adapted. Note, however, that a set-valued mapping $S$ can be adapted while its truncations $S_r$ are not. This happens e.g.\ when $S(\omega)=\{x\in\reals^n\mid \xi(\omega)x_0\le x_1\}$ with $\xi\in L^0\setminus L^\infty$.
\end{example}

\begin{proof}
Let $z\in\X$ such that $(z,\bar u)\in\dom Ef$ and let $r\ge\|z\|_{L^\infty}$. We have $z\in S_r$ and thus, $z^t\in S^t_r$ almost surely. When $S_r$ is adapted, Jensen's inequality, gives $E_tz^t\in S_r^t$ almost surely. Applying the measurable selection theorem to
\[
\omega\mapsto \{x\in\reals^n\mid x^t= [E_t z^t](\omega),\ x\in S_r(\omega)\},
\]
gives the existence of a $\bar z\in L^0$ such that $\bar z^t=E_t z^t$ and $\bar z\in S_r$ almost surely. The last condition now gives $(\bar z,\bar u)\in\dom Ef$.
\end{proof}

As noted in \cite{rw76}, any problem with essentially bounded domain can reduced to a problem satisfying relatively complete recourse. The following gives the details.

\begin{remark}[Induced constraints]
Assume that $S(\omega):=\dom f(\cdot,\bar u(\omega),\omega)$ is closed and essentially bounded and that $(x,\bar u)\in\dom Ef$ when $(x,\bar u)\in\dom f$ almost surely. The projections $S^t(\omega)$ of $S(\omega)$ are then closed and essentially bounded so, by \cite[Theorem~1.8]{tru91}, there exist $\F_t$-measurable mappings $D^t$ whose $\F_t$-measurable selections coincide with those of $S^t$. Clearly, every feasible $x\in\X_a$ has $x^t\in D^t$ almost surely. Thus, defining 
\[
\tilde f(x,u,\omega) := f(x,u,\omega) + \sum_{t=0}^T\delta_{D_t(\omega)}(x^t),
\]
we have $E\tilde f(x,\bar u)=Ef(x,\bar u)$ for every $x\in\X_a$. Moreover, $\tilde f$ satisfies the assumptions of \thref{ex:recinfty}. The random sets $D^t$ can be thought of as ``induced constraints'' arising form the requirement of feasible future recourse. If $S$ is adapted, we have $D^t=S^t$ and $f=\tilde f$.
\end{remark}



\section{Applications}\label{sec:app}

This section applies the general results above to specific instances of \eqref{spx}. The dual problems and optimality conditions were derived in \cite{pp22b}.

\subsection{Mathematical programming}

Consider the problem
\begin{equation}\label{mp}\tag{$MP$}
\begin{aligned}
&\minimize\quad & Ef_0(x)&\quad\ovr\ x\in\N,\\
  &\st\quad & f_j(x) &\le 0\quad j=1,\ldots,l\ a.s.,\\
   & & f_j(x) &= 0\quad j=l+1,\ldots,m\ a.s.
\end{aligned}
\end{equation}
where $f_j$ are convex normal integrands with $f_j$ affine for $j>l$. This extends the problem formulation from \cite{rw78} by relaxing the boundedness assumptions on the feasible set and by including the equality constraints.

Problem \eqref{mp} fits the general duality framework with $\bar u=0$ and
\[
f(x,u,\omega) = 
\begin{cases}
  f_0(x,\omega) & \text{if $x\in\dom H,\ H(x)+u\in K$},\\
  +\infty & \text{otherwise},
\end{cases}
\]
where $K=\reals_-^l\times\{0\}$ and $H$ is the $K$-convex random function defined by
\[
\dom H(\cdot,\omega)=\bigcap_{j=1}^m\dom f_j(\cdot,\omega)\quad\text{and}\quad H(x,\omega)=(f_i(x,\omega))_{j=1}^m.
\]
It was shown in \cite{pp22b} that if $\dom Ef\cap(\X\times\U)\ne\emptyset$, then the dual problem can be written as
\begin{equation}\label{dmp}\tag{$D_{MP}$}
\begin{aligned}
  &\maximize\ E\inf_{x\in\reals^n}\{f_0(x)+y\cdot H(x)-x\cdot p\}\ \ovr (p,y)\in\X_a^\perp\times\Y\\
  &\st\qquad\qquad y\in K^*\quad a.s.
\end{aligned}
\end{equation}
To get more explicit expressions for $f^*$ and the dual problem, additional structure is needed; see \cite{pp22b}.

In order to state our assumptions, we first write the problem in the form
\begin{equation*}
\begin{aligned}
  &\minimize\quad &  &Ef_0(x)\quad\ovr\ x\in\N,\\
  &\st &  &F(x)\in \reals^l_-\quad \text{a.s.},\\
  & &     &Ax=b\quad \text{a.s.},
\end{aligned}
\end{equation*}
where $F$ is an $\reals^l_-$-convex random function and $A$ is a random $n\times(m-l)$-matrix. Without loss of generality, we may assume that $b=(b_t)_{t=0}^T$ and 
\[
A=
\begin{bmatrix}
  A_{0,0} & 0 & \cdots & 0 \\
\vdots   & & & \vdots \\
 & & & 0\\
  A_{T,0} & & \cdots & A_{T,T}  
\end{bmatrix},
\]
where $b_t\in\reals^{m^2_t}$ and $A_{s,t}:\reals^{{n_s}}\to\reals^{m^2_t}$ with $m^2_1+\dots+m^2_T=m-l$. Similarly, we may assume that $F(x,\omega)=(F^0(x^0,\omega),\dots,F^T(x^T,\omega))$, where $F^t$ is a $\reals^{m^1_t}_-$-convex random function on $\reals^{n^t}$ with $m^1_1+\dots + m^1_T=l$.

Accordingly, we can write the normal integrand $f$ as
\begin{align*}
f(x,u,\omega) &=
\begin{cases}
  f_0(x,\omega) & \text{if $x\in\dom H,\ H(x)+u\in K$},\\
  +\infty & \text{otherwise},
\end{cases}\\
&=\begin{cases}
f_0(x,\omega)\quad &\text{if }F(x,\omega)+(u_j)_{j=1}^l\in\reals^l_-,\ A(\omega)x+(u_j)_{j=l+1}^m=b(\omega),\\
+\infty &\text{otherwise}.
\end{cases}
\end{align*}
We assume that $\X=L^\infty$ and that  $\U=L^\infty(\reals^l)\times \U^e$, where $\U^e$ is a space of $\reals^{m-l}$-valued random variables satisfying the assumptions in Section~\ref{sec:sdU}. We denote $A^t:=(A_{t,0},\ldots,A_{t,t})$. Since $\U^e$ is solid, we can write it accordingly as $\U^e=\U^e_0\times\cdots\times\U^e_T$, where $\U^e_t$ is a decomposable solid Fr\'echet space of $\reals^{m^2_t}$-valued random variables.

\begin{assumption}\thlabel{deMPass}\mbox{}
\begin{enumerate}[label=\Alph*]
\item
$Ef_0$ is finite on $\X$ and there exist $\bar x\in\X_a$ and $\epsilon>0$ such that
\begin{align*}
  f_j(\bar x+x) + \epsilon&\le 0\quad  \forall\ x \in\uball_\epsilon,\quad j=1,\ldots,l\ \text{a.s.} \\
  A(\bar x) &= b \quad \text{a.s.} 
\end{align*}
\item
  $A\X\subseteq\U^e$ and $A\X_a$ and $A\X$ are strongly closed in $\U^e$,
\item For every $t$, $A^t$, $b_t$ and $F^t$ are $\F_t$-measurable, and
\[
\{x_t \in \X_{t} \mid F^t(x^{t-1},x_t)\in\reals^{m^1_t}_-,\ A^t(x^{t-1},x_t) =b_t\ \text{a.s.}\}
\]
is nonempty for every $\F_{t-1}$-measurable $x^{t-1}\in \X^{t-1}$ such that $F^{t'}(x^{t'})\in\reals^{m^1_{t'}}_-$ and $A^{t'} x^{t'}=b_{t'}$ for all $t'<t$.
\end{enumerate}
\end{assumption}

The following extends the main results of \cite{rw78} relaxing the compactness assumptions made there and by allowing for affine equality constraints. The strategy of \cite{rw78} was to first relax the nonanticipativity constraint by using a shadow price of information $p$ and then to construct a dual variable $y$ via measurable selection arguments. The following employs the general theory of Sections~\ref{sec:sdp} and \ref{sec:recourse} to establish the existence of a dual optimal pair $(p,y)$ directly.

\begin{theorem}
Under \thref{deMPass}, $\inf\eqref{mp}=\sup\eqref{dmp}$ and the optimum in \eqref{dmp} is attained. In particular, a feasible $x\in\N$ solves \eqref{mp} if and only if there exists $(p,y)\in\N^\perp\times\Y$ feasible in \eqref{dmp} such that
\begin{equation*}
  \begin{gathered}
    p\in\partial_x[f_0+y\cdot H](x),\\
  H(x)\in K,\quad y\in K^*,\quad y\cdot H(x)=0
  \end{gathered}
\end{equation*}
almost surely.
\end{theorem}

\begin{proof}
Recall that, in \eqref{mp}, $\bar u=0$. It was shown in \cite[Section~6.1]{pp22b} that $(p,y)\in\partial f(x,\bar u)$ can be written as the scenario-wise optimality conditions given here. Thus, by \thref{thm:dexistence}, it suffices to show that  $\pos\dom\varphi$ is linear and closed  and \thref{ass:rec} holds. By \thref{deMPass}.A,
\begin{multline*}
  \dom Ef-(\bar x,0)=\{(z,u)\in\X\times\U\mid f_j(\bar x+z)+u_j\le 0\ j=1,\ldots,l\ a.s.\}\\
  \cap\{(z,u)\in\X\times\U\mid Az+(u_j)_{j=l+1}^m=0\},
\end{multline*}
where the first set on the right has $(0,0)$ in its strong interior. Thus, by \thref{lem:hulls},
\begin{align*}
  \pos(\dom Ef-(\bar x,0))  &= \{(x,u)\in\X\times\U\mid Ax+(u_j)_{j=l+1}^m=0\}.
\end{align*}
By \thref{lem:A}, the first property in $B$ implies that $A:\X\to\U$ is continuous so the problem is strictly feasible. By \thref{lem:posdomphi},
\begin{align*}
\pos\dom \varphi &= \pos(\dom Ef-(\bar x,0))-\X_a\times\{0\}\\
&= \{(z,u)\in\X\times\U\mid Az+(u_j)_{j=l+1}^m\in A\X_a\},
\end{align*} 
so $\pos\dom\varphi$ linear closed by B. Indeed, by \thref{lem:A}, the first property in $B$ implies that $A:\X\to\U$ is continuous.

The projection of $\pos(\dom Ef-(0,0)$ to $\U$ equals $L^\infty(\reals^l)\times A\X$, which is strongly closed by the last property in B. Thus, by \thref{lem:rec0}, \thref{ass:rec}.1 holds. As to \thref{ass:rec}.2, let $z\in\X$ with $(z,0)\in\dom Ef$. By \thref{deMPass}.C, $A^t$, $b_t$ and $F^t$ are  $\F_t$-measurable, so, by Jensen's inequality,
\begin{align*}
F^{t'}(E_tz^{t'})\in\reals_-^{t'} \quad\forall t'\le t,\\
A^{t'}(E_tz^{t'})=b_{t'}\quad\forall t'\le t.
\end{align*}
Under \thref{deMPass}.C, $E_tz^{t}$ can be extended to a $\bar z$ such that $\bar z^t=E_tz^t$ and $(\bar z,0)\in\dom Ef$.
\end{proof}

The following gives sufficient conditions for the closedness conditions in \thref{deMPass}.B. Recall that, if $A$ is a matrix with full row rank, its Moore--Penrose inverse is given by
\[
A^\dagger = A^*(A^*A)^{-1}.
\]

\begin{remark}
Assume that $A\X\subseteq\U^e$ and that, for every $t$, $A^t$ is $\F_t$-measurable and
\begin{equation}\label{eq:Att}
\U^e_t=A_{t,t}\X_t.
\end{equation}  
Then $B$ holds. If there are no inequality constraints and $b_t$ are $\F_t$-measurable, then $C$ holds.

Condition \eqref{eq:Att} holds trivially if $\U^e_t$ is the range space of the mapping $\A_{t,t}:\X_t\to L^0$ defined pointwise by $(\A_{t,t}x_t)(\omega)=A(\omega)x_t(\omega)$. By \thref{ex:scaled}, such a space $\U^e_t$ satisfies the assumptions of the present chapter (i.e.\ it is a solid Fr\'echet space and its dual is a direct sum of measurable functions and singular elements) provided $|A_{t,t}||A_{t,t}^\dagger|\in L^\infty$ and
\[
\X_t=\{x_t\in L^0\mid |x_t|\in\X_t^0\}
\]
for solid decomposable spaces of scalar-valued random variables $\X^0$.
\end{remark}

\begin{proof}
To verify B, it suffices to prove that $A\X_a=\U^e_a$ and $A\X=\U^e$. Here $\U^e_a$ is the set of $u\in\U^e$ such that the component $u_t$ of $u$ corresponding to $A^t$ is $\F_t$-measurable. The block-triangular structure and the $\F_t$-measurability of $A^t$ imply $A\X_a\subset \U^e_a$. Given $u\in\U^e_a$, assume that there exists $x^{t-1}\in\X^{t-1}_a$ such that $A^{t'}x^{t'}+u_{t'}=0$ for all $t'\le t-1$. This holds trivially for $t=0$. Since $A\X\subseteq\U^e$, condition \eqref{eq:Att} gives the existence of an $x_t\in \X_t$ with $A^tx^t+u_t=0$. Taking conditional expectations, $A^tE_t x^t+u_t=0$, so there exists $x^t \in \X^t_a$ with $A^{t'}x^{t'}+u_{t'}=0$ for all $t'\le t$. By induction, $A\X_a=\U^e_a$. The equality, $A\X=\U^e$ is proved similarly. Assumption A implies $b\in\U^e$ so, in the absence of inequality constraints, the above argument also shows that $C$ holds.
\end{proof}





\subsection{Optimal stopping}

Let $R$ be a real-valued adapted stochastic process with $R_t\in L^1$ for all $t$ and consider the {\em optimal stopping problem}
\begin{equation}\label{os}\tag{$OS$}
  \maximize\quad ER_\tau\quad\ovr \tau\in\T,
\end{equation}
where $\T$ is the set of {\em stopping times}, i.e.\ measurable functions $\tau:\Omega\to\{0,\ldots,T+1\}$ such that $\{\omega\in\Omega\mid \tau(\omega)\le t\}\in\F_t$ for each $t=0,\ldots,T$. Choosing $\tau=T+1$ is interpreted as not stopping at all. The problem
\begin{equation}\label{ros}\tag{$ROS$}
\begin{aligned}
&\maximize\quad & & E\sum_{t=0}^TR_t x_t\quad\ovr x\in\N,\\
&\st\quad & & x\ge 0,\ \sum_{t=0}^Tx_t\le 1\quad \text{a.s.}
\end{aligned}
\end{equation}
is the convex relaxation of \eqref{os} in sense that their optimal values coincide and the extreme points of the feasible set of \eqref{ros} can be identified with $\T$; see \cite[Section~5.2]{pp22}. 

Problem \eqref{ros} fits the general duality framework with $n_t=1$, $\X=L^\infty$, $\V=L^1$, $m=1$, $\U=L^\infty$, $\Y=L^1$,
\[
f(x,u,\omega) =
\begin{cases}
-\sum_{t=0}^Tx_tR_t(\omega) & \text{if $x\ge 0$ and $\sum_{t=0}^Tx_t+u\le 0$},\\
+\infty & \text{otherwise}
\end{cases}
\]
and $\bar u=-1$. It was shown in \cite[Section~6.2]{pp22b} that the dual of \eqref{ros} can be written as
\begin{equation}\label{dos}\tag{$D_{OS}$}
\begin{aligned}
&\minimize\quad Ey\quad\ovr(p,y)\in\X_a^\perp\times \Y_+\\
&\st\quad p_t+R_t\le y\quad t=0,\ldots,T\ \text{a.s.}
\end{aligned}
\end{equation}
As noted in \cite[Section~6.2]{pp22b}, a pair $(p,y)\in \V\times \Y$ solves \eqref{dos} if and only if $p_t=y-E_t y$ and the process $y_t:=E_ty$ solves the ``reduced dual''
\begin{align*}
&\minimize\quad Ey_0\quad\ovr y\in\M^\Y_+\\
&\st\quad R_t\le y_t\quad t=0,\ldots,T\ \text{a.s.},
\end{align*}
where $\M^\Y_+$ is the cone of nonnegative martingales $y$ with $y_T\in\Y$.

\begin{theorem}
We have $\inf\eqref{os}=\sup\eqref{dos}$ and the optimum in \eqref{dos} is attained. In particular, a stopping time $\tau\in\T$ is optimal if and only if there exists $(p,y)\in\X_a^\perp\times\Y$ such that $p_t+R_t\le y$ for all $t$ and $p_\tau+R_\tau=y$ almost surely. This is equivalent to the existence of a martingale $y$ such that $R_t\le y_t$ for all $t$ and $R_\tau=y_\tau$ almost surely.
\end{theorem}

\begin{proof}
Recall that in \eqref{ros}, $\bar u=-1$. By \cite[Theorem~38 and Example~39]{pp22b}, it suffices to show that $\varphi$ is subdifferentiable at $(0,\bar u)$. Thus, by \thref{thm:dexistence},  it suffices to show that $\pos(\dom\varphi-(0,\bar u))$ is linear and closed and \thref{ass:rec} holds. We have
\[
\dom Ef = \{(x,u)\in L^\infty\times L^\infty\mid x\ge 0,\ \sum_{t=0}^Tx_t+u\le 0\}
\]
and $\dom\varphi = L^\infty_-$. Thus $\pos(\dom\varphi-(0,\bar u))$ is the whole space so it is linear and closed. We have $\dom_u Ef=L^\infty_-$ so \thref{ass:rec}.1 holds. The mapping $S(\omega):=\dom f(\cdot,\bar u(\omega),\omega)=\{x\in\reals^n\mid x\ge 0,\ \sum_{t=0}^Tx_t\le 1\}$ is deterministic and, in particular, adapted so \thref{ass:rec}.2 holds by \thref{ex:recinfty}.
\end{proof}

\subsection{Optimal control}

Consider the optimal control problem
\begin{equation}\label{oc}\tag{$OC$}
\begin{aligned}
&\minimize\quad & & E\left[\sum_{t=0}^{T} L_t(X_t,U_t)\right]\quad\ovr\ (X,U)\in\N,\\
&\st\quad & & \Delta X_{t}=A_t X_{t-1} +B_t U_{t-1}+W_t\quad t=1,\dots,T
\end{aligned}
\end{equation}
where the {\em state} $X$ and the {\em control} $U$ are processes with values in $\reals^N$ and $\reals^M$, respectively, $A_t$ and $B_t$ are $\F_t$-measurable random matrices, $W_t$ is an $\F_t$-measurable random vector and the functions $L_t$ are convex normal integrands. The linear constrains in \eqref{oc} are called the {\em system equations}.

The problem fits the general duality framework with $x=(X,U)$, $\bar u=(W_t)_{t=1}^T$ and
\[
f(x,u,\omega)=\sum_{t=0}^{T}L_t(X_t,U_t,\omega) + \sum_{t=1}^{T}\delta_{\{0\}}(\Delta X_{t}- A_t(\omega)X_{t-1}-B_t(\omega)U_{t-1}-u_t).
\]
As in \cite[Section~6.3]{pp22b}, we assume that
\[
\begin{aligned}
  \X_t&=\S\times\C,&\quad \U_t&=\S\\
  \V_t&=\S'\times\C', &\quad \Y_t&=\S'
\end{aligned}
\]
where $\S$ can $\C$ are solid decomposable spaces in separating duality with $\S'$ and $\C'$, respectively. In order to apply the general theory of Sections~\ref{sec:sdp} and \ref{sec:recourse}, we assume, in addition, that $E_t\S\subset\S$, $E_t\C\subseteq\C$ for all $t$ and that $\S$ and $\C$ are endowed with Fr\'echet topologies under which their topological duals can be expressed as
\[
\S^*=\S'\oplus(\S')^s\quad\text{and}\quad\C^*=\C'\oplus(\C')^s,
\]
where $(\S')^s$ and $(\C')^s$ are singular elements of $\S^*$ and $\C^*$, respectively. These conditions simply mean that the spaces $\X$ and $\U$ satisfy the assumptions made in Sections~\ref{sec:sdp} and \ref{sec:recourse}.

It was shown in \cite[Section~6.3]{pp22b} that if $\dom Ef\cap(\X\times\U)\ne\emptyset$, then the dual of \eqref{oc} can be written as
\begin{equation}\label{doc}\tag{$D_{OC}$}
\begin{aligned}
&\maximize & & E\left[\sum_{t=1}^{T}W_t\cdot y_t-\sum_{t=0}^{T} L^*_t(p_t - (\Delta y_{t+1}+A^*_{t+1}y_{t+1},B^*_{t+1}y_{t+1}))\right]\\
&\ovr & & (p,y)\in\X_a^\perp\times\Y.
\end{aligned}
\end{equation}

\begin{assumption}\thlabel{deOCass}
For $t=0,\ldots,T$
\begin{enumerate}[label=\Alph*]
\item\label{OCassB} $A_{t}\S\subseteq\S$ and $B_{t}\C\subseteq\S$,
\item\label{OCassC}  $EL_t$ is proper on $\S\times\C$, and there exists a feasible $(\bar X,\bar U)\in\X_a$ such that
\[
\K_t:=\pos(\dom EL_t-(\bar X_t,\bar U_t))\quad t=0,\ldots,T
\]
are linear and strongly closed.
\item\label{OCassDpos}
\begin{enumerate}
\item $E_{t}\K_{t'}\subset \K_{t'}$ for all $t'\le t$
\item for every $X_t\in\S$, there exists $U_t\in\C$ with $(X_t,U_t)\in\K_t$. 
\end{enumerate}
\end{enumerate}
\end{assumption}

Clearly, \thref{deOCass}.\ref{OCassC} and \ref{OCassDpos} hold if the functions $EL_t$ are finite on $\S\times\C$. \thref{deOCass}.\ref{OCassB} was used also in \cite[Section~6.3]{pp22b} to derive a reduced dual problem. Part (a) of C holds automatically if each $L_t$ is $\F_t$-measurable.

Let $\Y_a$ be the adapted elements of $\Y$. By \thref{lem:EGcont}, $\ap y\in\Y_a$ for every $y\in\Y$, where  $(\ap y)_t:= E_t y_t$ is the {\em adapted projection} of $y$.  As noted in \cite[Section~6.3]{pp22b}, if \thref{deOCass}.A holds, each $L_t$ is $\F_t$-measurable and $EL_t$ is proper on $\S\times\C$, then a pair $(p,y)$ solves \eqref{doc} if and only if $\ap y$ solves the reduced dual
\begin{equation*}
\begin{aligned}
&\maximize\quad & & E\left[\sum_{t=1}^T W_t\cdot y_t -\sum_{t=0}^{T} [L^*_t(-E_t(\Delta y_{t+1}+A^*_{t+1}y_{t+1}, E_tB^*_{t+1}y_{t+1}))]\right]\quad\ovr\  y\in\Y_a
\end{aligned}
\end{equation*}
and 
\[
p_t = (\Delta y_{t+1}+A^*_{t+1}y_{t+1},B^*_{t+1}y_{t+1})- E_t(\Delta y_{t+1}+A^*_{t+1}y_{t+1},B^*_{t+1}y_{t+1}).
\]

\begin{lemma}\thlabel{lem:ocrcore}
Under \thref{deOCass}, $\pos(\dom\varphi-(0,\bar u))$ is linear and closed and \thref{ass:rec}.1 holds.
\end{lemma}

\begin{proof}
This fits the format of \thref{ex:rcoregA} with $f(x,u) = g(x,Ax+u)$, where $g(x,u):=\sum L_t(x_t)+\delta_{\{0\}}(u)$, $A(x)_t=-\Delta X_{t}+ A_tX_{t-1}+B_t U_{t-1}$ and $\pi'(x,u):=(\ap x,\ap u)$. By \thref{deOCass}.\ref{OCassC}, 
\[
\K:=\pos(\dom Eg-(\bar x,A\bar x+\bar u))=\{(x,u)\in\X\times\U \mid x_t\in\K_t\ \forall t,\, u=0\}
\]
is linear and strongly closed. By \thref{deOCass}.C, $\pi' \K\subseteq\K$. By \thref{thm:ssd,ex:rcoregA}, it thus suffices to show that 
\[
(I,A)\X_a\subset\rge\pi'\subseteq\K+(I,A)\X_a.
\]
We have $\rge\pi'=\X_a\times\U_a$. Since $A_t$ and $B_t$ are $\F_t$-measurable, the first inclusion is clear. As to the second, let $(x,u)\in\X_a\times\U_a$. We have $(x,u)=(I,A)x+ (0,u-Ax)$, where $Ax \in\U_a$, so it suffices to show that $(0,u)\in\K+(I,A)\X_a$, i.e., there exists $x\in\X_a$ with $(X_t,U_t)\in\K_{t}$ and $\Delta X_{t}= A_{t} X_{t-1}+B_{t} U_{t-1}+u_{t}$ for all $t$.

Assume that $(X^{t-1},U^{t-1})\in\X^{t-1}_a$ with 
\[
\Delta X_{t'}= A_{t'} X_{t'-1}+B_{t'} U_{t'-1}+u_{t'}\quad \forall t'=0\dots t-1
\]
and $(X_{t-1},U_{t-1})\in\K_{t-1}$. By \thref{deOCass}.A, the $\F_t$-measurable $X_t$ defined by
\[
\Delta X_{t}:= A_{t} X_{t-1}+B_{t} U_{t-1}+u_{t}
\]
belongs to $\S$. By \thref{deOCass}.\ref{OCassDpos}(b), there exists $\tilde U_t$ such that $(X_t,\tilde U_t)\in\K_t$.  Since $X_t$ is $\F_t$-measurable, \thref{deOCass}.\ref{OCassDpos}.(a) gives $(X_t,U_t)\in\K_t$, where $U_t:=E_t\tilde U_t$. Thus the claim follows by induction.

To verify \thref{ass:rec}.1, we apply \thref{ex:rec1}. It suffices to show that $\{u\in\U\mid (0,u)\in\K-\rge A\}=\U$. This means that, for every $u\in\U$, there exists $x\in\X$ with $(X_t,U_t)\in\K_{t}$ and $\Delta X_{t}= A_{t} X_{t-1}+B_{t} U_{t-1}+u_{t}$ for all $t$. The existence of such $x$ follows by repeating the arguments in the above paragraph.
\end{proof}

\begin{assumption}\thlabel{OCassD}
For all $t$,
\begin{enumerate}
\item $E_{t}\dom EL_{t'}\subset\dom EL_{t'}$ for all $t'\le t$
\item for every $X_t\in\S$, there exists $U_t\in\C$ with $(X_t,U_t)\in\dom EL_t$.\end{enumerate}
\end{assumption}

\thref{OCassD} implies \thref{deOCass}.C. Note also that \thref{OCassD}.1 holds automatically if each $L_t$ is $\F_t$-measurable.

\begin{theorem}
Under \thref{deOCass,OCassD}, problem \eqref{oc} is feasible for any $W\in\U_a$, $\inf\eqref{oc}=\sup\eqref{doc}$ and the optimum in \eqref{doc} is attained. In particular, a feasible $(X,U)$ is optimal in \eqref{oc} if and only if there exists a dual feasible $(p,y)\in\N^\perp\times\Y$ such that
\begin{align*}
p_t-(\Delta y_{t+1}+A^*_{t+1}y_{t+1},B^*_{t+1} y_{t+1})\in\partial L_t(X_t,U_t)
\end{align*}
almost surely. If each $L_t$ is $\F_t$-measurable, this is equivalent to the existence of a $y\in\Y_a$ feasible in the reduced dual such that 
\begin{align*}
-E_t(\Delta y_{t+1}+A^*_{t+1}y_{t+1},B^*_{t+1} y_{t+1})\in\partial L_t(X_t,U_t)
\end{align*}
almost surely.
\end{theorem}

\begin{proof}
Recall that in \eqref{oc}, $\bar u=W$.  By \cite[Theorem~41 and Remark~45]{pp22b}, it suffices to show that $\varphi$ is subdifferentiable at $(0,\bar u)$.  Thus, by \thref{thm:dexistence,lem:ocrcore}, it suffices to show that \thref{ass:rec}.2 holds and that $(0,\bar u)\in\dom\varphi$ for all $\bar u\in\U_a$. We start with the latter.
  
Let $\bar u\in\U_a$ and assume that $(X^{t-1},U^{t-1})\in\X^{t-1}_a$ with 
\[
\Delta X_{t'}= A_{t'} X_{t'-1}+B_{t'} U_{t'-1}+\bar u_{t'}\quad \forall t'=0\dots t-1
\]
and $(X_{t-1},U_{t-1})\in\dom EL_{t-1}$. By \thref{deOCass}, the $\F_t$-measurable $X_t$ defined by
\[
\Delta X_{t}= A_{t} X_{t-1}+B_{t} U_{t-1}+\bar u_{t}
\]
belongs to $\S$. By \thref{OCassD}.2, there exists $\tilde U_t$ such that $(X_t,\tilde U_t)\in\dom EL_t$. Since $X_t$ is $\F_t$-measurable, \thref{OCassD}.1 gives $(X_t,U_t)\in\dom EL_t$, where $U_t:=E_t\tilde U_t$. By induction, there exists $x\in\X_a$ such that $(x,\bar u)\in\dom Ef$, so $(0,\bar u)\in\dom\varphi$. 

To verify \thref{ass:rec}.2,  let $(z,\bar u)\in\dom Ef$. For all $t'\le t$,  $E_t z_{t'}\in \dom EL_{t'}$, by \thref{OCassD}.1, and, by \cite[Lemma~44]{pp22b},
\begin{align*}
E_{t} X_{t'} &= E_t X_{t'-1} +A_{t'} E_{t} X_{t'-1}+B_{t'}  E_{t} U_{t'-1}+ \bar u_{t'}.
\end{align*}
An induction argument similar to the above paragraph gives  $\bar z\in\X$ with $(\bar z,\bar u)\in\dom Ef$ and $\bar z^t=E_t z^t$, so \thref{ass:rec}.2 holds.
\end{proof}

\begin{example}[Bounded strategies]\thlabel{OCbdd}
Assume that $\S=L^\infty$ and $\C=L^\infty$ are endowed with the usual norm-topologies. If $EL_t$ is finite on
\[
\D_t := \{(X,U)\in L^\infty\times L^\infty\mid\exists\epsilon>0:\ (X_t,U_t)+\epsilon\uball\subset\dom L_t\quad P\text{-a.s.}\}
\]
and if there is a feasible point $(\bar X,\bar U)$ such that $(\bar X_t,\bar U_t)\in\D_t$ for all $t$, then, by \thref{thm:rinte}, \thref{deOCass}.\ref{OCassC} holds.
\end{example}

\subsection{Problems of Lagrange}

Consider the problem
\begin{equation}\label{lagrange}\tag{$L$}
\minimize\quad E\sum_{t=0}^T K_t(x_t,\Delta x_t)\quad\ovr x\in\N,
\end{equation}
where $x$ is a process of fixed dimension $d$, $K_t$ are convex normal integrands and $x_{-1}:=0$. Problem~\eqref{lagrange} can be thought of as a discrete-time version of a problem studied in calculus of variations. Other problem formulations have $K_t(x_{t-1},\Delta x_t)$ instead of $K_t(x_t,\Delta x_t)$ in the objective, or an additional term of the form $Ek(x_0,x_T)$ (see \cite{rw83}), all of which fit the general format of stochastic optimization.

This fits the general duality framework with $\bar u=0$ and
\[
f(x,u,\omega) = \sum_{t=0}^T K_t(x_{t},\Delta x_t+u_t,\omega).
\]
As in \cite[Section~6.4]{pp22b}, we assume that
\[
\X_t=\S,\quad \V_t=\S',\quad \U=\X,\quad\Y=\V,
\]
where $\S$ is solid decomposable space in separating duality with $\S'$. In addition, we assume that $E_t\S\subset\S$ for all $t$ and that $\S$ is endowed with a Fr\'echet topology under which its topological dual can be expressed as
\[
\S^*=\S'\oplus(\S')^s,
\]
where $(\S')^s$ are singular elements of $\S^*$. These conditions simply mean that the spaces $\X$ and $\U$ satisfy the assumptions made in Sections~\ref{sec:sdp} and \ref{sec:recourse}.

It was shown in \cite[Section~6.4]{pp22b} that if $\dom Ef\cap(\X\times\U)\ne\emptyset$, then the dual problem can be written as
\begin{equation}\label{dl}\tag{$D_L$}
\maximize\quad\quad E [-\sum_{t=0}^T K_t^*(p_{t}+\Delta y_{t+1},y_t)]\quad\ovr y\in\Y,p\in \X_a^\perp
\end{equation}
where $y_{T+1} := 0$.

As noted in \cite[Section~6.4]{pp22b}, if each $K_t$ is $\F_t$-measurable and $EK_t$ is proper on $\S\times\S$, then a pair $(p,y)$ solves the dual problem if and only if $\ap y$ solves the reduced dual problem
\begin{equation*}
\begin{aligned}
&\maximize\quad\quad  & & E [-\sum_{t=0}^T K_t^*(E_t\Delta y_{t+1},y_t)]\quad\ovr y\in\Y_a.
\end{aligned}
\end{equation*}
and $p_t=E_ty_{t+1}-y_{t+1}$.

\begin{assumption}\thlabel{deBolzaAss} \mbox{}
\begin{enumerate}[label=\Alph*]
\item\label{BolzaAssB} each $EK_t$ is proper lsc on $\S\times\S$.
\item\label{BolzaAssC} there exists $\bar x\in\X_a$ such that, for all $t$, 
\[
\K_t:=\pos(\dom EK_t-(\bar x_t,\Delta \bar x_t))
\]
is linear and strongly closed,
\item\label{BolzaAssDpos} for every  $t$, 
\begin{enumerate}
\item $E_t \K_{t'}\subset\K_{t'}$ for all $t'\le t$,
\item for every $u_t,x_{t-1}\in\S$, there exists $x_t\in\S$ such that $(x_t,\Delta x_t+u_t)\in \K_t$.
\end{enumerate}
\end{enumerate}
\end{assumption}

\begin{lemma}\thlabel{lagrangercore}
Under \thref{deBolzaAss}, $\pos\dom\varphi$ is linear and closed and \thref{ass:rec}.1 holds.
\end{lemma}

\begin{proof}
This fits the format of \thref{ex:rcoregA} with $f(x,u) = g(x,Ax+u)$, where $g(x,u):=\sum_{t=0}^T K_t(x_t,u_t)$, $A(x)_t :=\Delta x_t$ and $\pi'(x,u):=(\ap x,\ap u)$. We have
\[
\K := \pos(Eg-(\bar x,A\bar x+\bar u))=\prod_{t=0}^T\K_t,
\]
which is linear and strongly closed, by \thref{deBolzaAss}.B. \thref{deBolzaAss}.C(a) implies $\pi'\K\subseteq\K$. By \thref{thm:ssd,ex:rcoregA}, it thus suffices to show that 
\[
(I,A)\X_a\subseteq\rge\pi'\subseteq\K-(I,A)\X_a.
\]
The first inclusion is clear. 

The second inclusion means that, for every $(x,u)\in\rge\pi'=\X_a\times\X_a$, there exists $\tilde x\in\X_a$ such that $(x_t+\tilde x_t,u_t+\Delta \tilde x_t)\in \K_t$ for all $t$. By change of variables, this means for every $u\in\X_a$, there exists $\tilde x\in\X_a$ such that $(\tilde x_t,u_t+\Delta \tilde x_t)\in \K_t$ for all $t$. Assume that there is $\tilde x^t\in\X^t_a$ such that  $(x_{t'},u_{t'}+\Delta \tilde x_{t'})\in \K_t$ for all $t'\le t$. By \thref{deBolzaAss}.C(b), there exists $\tilde x_{t+1}\in\S$ such that  $(\tilde x_{t+1},u_{t+1}+\Delta \tilde x_{t+1})\in \K_{t+1}$. By \thref{deBolzaAss}.C(a) and $\F_{t+1}$ measurability of $x_{t+1}$ and $u_{t+1}$, we can choose $\tilde x_{t+1}$ as  $\F_{t+1}$-measurable. Thus the second inclusion follows by induction on $t$.


To verify \thref{ass:rec}.1, we apply \thref{ex:rec1}. It suffices to show that $\{u\in\U\mid (0,u)\in\K-\rge A\}=\U$. This means that, for every $u\in\U$, there exists $x\in\X$ with $(x_t,\Delta x_t+u_t)\in\K_{t}$ for all $t$. The existence of such $x$ follows by repeating the induction argument in the above paragraph.
\end{proof}

\begin{assumption} \thlabel{BolzaAssD} 
For every $t$, 
\begin{enumerate}
\item $E_t \dom EK_{t'}\subset \dom EK_{t'}$ for all $t'\le t$,
\item for every $\F_t$-measurable $x^t\in\X^t$ such that $(x_t,\Delta x_t)\in\dom EK_t$, there exists $x_{t+1}\in\X_{t+1}$ with $(x_{t+1},\Delta x_{t+1})\in\dom EK_{t+1}$.
\end{enumerate}
\end{assumption}

When
\[
\dom EK_t=\{(x_t,u_t)\in \S\times\S\mid (x_t,u_t)\in\dom K_t\ a.s.\},
\]
1 holds, by Jensen, as soon as $\dom K_t$ is $\F_t$-measurable. Moreover, 2 holds if $\dom_1 EK_t$ or $\dom_2 EK_t$ is the whole space.

Like in \cite[Section~6.4]{pp22b}, we will write the optimality conditions in terms of the associated {\em Hamiltonians}
\[
H_t(x_{t},y_t,\omega) := \inf_{u_t\in\reals^d}\{K_t(x_{t},u_t,\omega)-u_t\cdot y_t\}.
\]

\begin{theorem}
Under \thref{deBolzaAss,BolzaAssD}, $\inf\eqref{lagrange}=\sup\eqref{doc}$ and the optimum in \eqref{dl} is attained. In particular, a feasible $x$ is optimal in \eqref{lagrange} if and only if there exists a dual feasible $(p,y)\in\N^\perp\times\Y$ such that
    \begin{align*}
      \begin{split}
        p_{t} + \Delta y_{t+1}&\in \partial_{x} H_t(x_{t},y_t),\\
        \Delta x_t &\in \partial_{y}[-H_t](x_{t},y_t),
      \end{split}
    \end{align*}
almost surely. If, in addition, each $K_t$ is $\F_t$ measurable, this is equivalent to the existence of a $y\in\Y_a$ feasible in the reduced dual such that 
\begin{align*}
\begin{split}
E_t\Delta y_{t+1}&\in \partial_{x} H_t(x_{t},y_t),\\
\Delta x_t &\in \partial_{y}[-H_t](x_{t},y_t),
\end{split}
\end{align*}
almost surely.
\end{theorem}

\begin{proof}
Recall that in \eqref{lagrange}, $\bar u=0$. By \cite[ Theorem~49 and Remark~51]{pp22b}, it suffices to show that $\varphi$ is subdifferentiable at the origin. Thus, by \thref{thm:dexistence,lagrangercore}, it suffices to show that \thref{ass:rec}.2 holds.

Let $z\in\X$ be such that $Ef(z,0)<\infty$. Let $\bar z^t=E_tz^t$. By \thref{BolzaAssD}.1, $(\bar z_{t'},\Delta\bar z_{t'})=E_t(z_{t'},\Delta z_{t'})\in\dom EK_{t'}$ for all $t'\le t$. By \thref{BolzaAssD}.2, there is an $\tilde z_{t+1}\in\X_{t+1}$ with $(\tilde z_{t+1},\Delta\tilde z_{t+1})\in\dom EK_{t+1}$. Let $\bar z_{t+1}=E_{t+1}\tilde z_{t+1}$. By \thref{BolzaAssD}.1, $(\bar z_{t+1},\Delta\bar z_{t+1})\in\dom EK_{t+1}$. Repeating until $T$, we find a $\bar z\in\X$ that satisfies the conditions in \thref{ass:rec}.2.
\end{proof}

\begin{example}[Bounded strategies]\thlabel{ex:scbdd}
Assume that $\S=L^\infty$ is endowed with the usual norm-topology. If $EK_t$ is finite on
\[
\D_t := \{(x_t,u_t)\in L^\infty\times L^\infty\mid\exists\epsilon>0:\ (x_t,u_t)+\epsilon\uball\subset\dom K_t\quad P\text{-a.s.}\}
\]
and if there is a feasible $\bar x$ such that $(\bar x_t,\Delta \bar x_t)\in\D_t$ for all $t$, then, by \thref{thm:rinte}, \thref{deBolzaAss}.\ref{BolzaAssC} holds.
\end{example}

\subsection{Financial mathematics}\label{sec:fm5}

This section analyzes the semi-static hedging problem from \cite{pp22b}. Let $s=(s_t)_{t=0}^T$ be an adapted $\reals^J$-valued stochastic process describing the unit prices of a finite set $J$ of perfectly liquid tradeable assets. We also assume that there is a finite set $K$ of derivative assets that can be bought or sold at time $t=0$ and that provide random payments $C^j\in L^0$, $j\in K$ at time $t=T$. We denote $C=(C^j)_{j\in K}$. The cost of buying a derivative portfolio $\bar x\in\reals^{K}$ at the best available market prices is denoted by $S_{-1}(x_{-1})$. Such a function is convex and lsc; see e.g.\ \cite{pen11,pen11b}. 

Consider the problem of finding a dynamic trading strategy $x=(x_t)_{t=0}^T$ in the liquid assets $J$ and a static portfolio $x_{-1}$ in the derivatives $K$ so that their combined revenue provides the ``best hedge'' against the financial liability of delivering a random amount $c\in L^0$ of cash at time $T$. If we assume that cash (or another numeraire asset) is a perfectly liquid asset that can be lent and borrowed at zero interest rate, the problem can be written as
\begin{equation*}\label{ssh}\tag{$SSH$}
\begin{aligned}
&\minimize\quad & & EV\left(c - \sum_{t=0}^{T-1}x_t\cdot\Delta s_{t+1} - C\cdot x_{-1} + S_{-1}(x_{-1})\right)\ \ovr\ x\in\N,x_{-1}\in\reals^{K},\\
  &\st\quad & & x_t \in D_t\quad t=0,\ldots,T-1\ a.s.,
\end{aligned}
\end{equation*}
where $V$ is a random ``loss function'' on $\reals$ and $D_t$ is a random $\F_t$-measurable set describing possible portfolio constraints. More precisely, the function $V$ is a convex normal integrand such that $V(\cdot,\omega)$ nondecreasing and nonconstant for all $\omega$. We will assume $D_T=\{0\}$, which means that all positions have to be closed at the terminal date. 


As soon as $c\in\U$, problem~\eqref{ssh} fits the general duality framework with the time index running from $-1$ to $T-1$, $\F_{-1}=\{\Omega,\emptyset\}$, $\bar u=c$ and
\[
f(x,u,\omega)=V\left(u - \sum_{t=0}^{T-1}x_t\cdot\Delta s_{t+1}(\omega) - C(\omega)\cdot x_{-1} + S_{-1}(x_{-1}),\omega\right)+ \sum_{t=0}^{T-1}\delta_{D_t(\omega)}(x_t,\omega).
\]
We will assume that $S_{-1}(0)=0$ and $0\in D_t$ almost surely for all $t$. It was shown in \cite[Section~6.5]{pp22b} that as soon as $EV$ is proper on $\U$, then the dual of \eqref{ssh} can be written as
\begin{equation}\label{dssh}\tag{$D_{SSH}$}
  \begin{aligned}
    &\maximize_{p\in \X_a^\perp, y\in\Y}\quad & & E\left[cy - V^*(y) - \sum_{t=0}^{T-1}\sigma_{D_t}(p_t+y\Delta s_{t+1}) - (yS_{-1})^*(p_{-1}+yC)\right].
  \end{aligned}
\end{equation}

We assume that $\X=L^\infty$ and that  $\U$ is a space of real-valued random variables satisfying the assumptions in Section~\ref{sec:sdU}. 

\begin{assumption}\thlabel{ass:fmds}\mbox{ }
\begin{enumerate}[label=\Alph*]
\item $\Delta s^j_t,C_{\bar j}\in \U$ for all $j\in J, \bar j\in K$ and $t=0,\ldots,T$,
\item $EV$ is finite on $\U$, $EV^*$ is proper on $\Y$,
\item there exists a feasible $\bar x\in\X_a$ and an $\epsilon>0$ such that $\bar x_{-1}\in\inte\dom S_{-1}$ and $\uball_\epsilon(\bar x_t) \subseteq D_t$ almost surely,
\item $D_t$ is $\F_t$-measurable for $t=0,\ldots,T$ and $S_{-1}(z_{-1})\in\U$ for all $z_{-1}\in L^\infty(\dom S_{-1})$.
\end{enumerate}
\end{assumption}

As noted in \cite[Section~6.5]{pp22b}, if $E_t\U\subseteq\U$ and $\Delta s_{t+1}\in\U$, then  $(p,y)$ solves \eqref{dssh} if and only if $y$ solves the reduced dual problem
\begin{equation*}
  \begin{aligned}
    &\maximize_{y\in\Y}\quad & & E\left[cy - V^*(y) - \sum_{t=0}^{T-1}\sigma_{D_t}(E_t[y\Delta s_{t+1}]) - (E[y]S_{-1})^*(E[y\bar c])\right].
  \end{aligned}
\end{equation*}
and
\[
p_{-1}:=\frac{E[y\bar c]}{E[y]}y-y\bar c\quad\text{and}\quad p_t=E_t[y\Delta s_{t+1}]-y\Delta s_{t+1}\quad t=0,\ldots,T-1.
\]

\begin{theorem}\thlabel{thm:fm5}
Under \thref{ass:fmds}, $\inf\eqref{ssh}=\sup\eqref{dssh}$ and the optimum in \eqref{dssh} is attained. In particular, a feasible $x\in\N$ solves \eqref{ssh} if and only if there exists $(p,y)\in\N^\perp\times\Y$ feasible in \eqref{dssh} such that
\begin{align*}
y &\in\partial V(u-\sum_{t=0}^{T-1} x_t\Delta s_{t+1}-C\cdot x_{-1}+S_{-1}(x_{-1})),\\
p_t+y\Delta s_{t+1} &\in N_{D_t}(x_t)\quad t=0,\dots,T,\\
p_{-1}+yC&\in\partial(yS_{-1})(x_{-1})
\end{align*}
almost surely. This holds if and only if
\begin{align*}
y &\in\partial V(u-\sum_{t=0}^{T-1} x_t\Delta s_{t+1}-C\cdot x_{-1}+S_{-1}(x_{-1})),\\
E_t[y\Delta s_{t+1}] &\in N_{D_t}(x_t)\quad t=0,\dots,T,\\
\frac{E[yC]}{E[y]}y&\in\partial (yS_{-1})(x_{-1}),
\end{align*}
where the fraction is interpreted as $0$ if $E[y]=0$.
\end{theorem}

\begin{proof}
It was shown in \cite[Theorem~56 and Remark~58]{pp22b} that the scenario-wise optimality condition $(p,y)\in\partial f(x,\bar u)$ can be written as the scenario-wise optimality conditions given here.  Thus, by \thref{thm:dexistence}, it suffices to show that $\pos(\dom\varphi-(0,\bar u))$ is linear and closed and \thref{ass:rec} holds.

We have
\[
\varphi(z,u)\le EV(u-\sum_{t=0}^T ((\bar x_t+z_t)\cdot\Delta s_{t+1})-C\cdot(\bar x_{-1}+z_{-1})+S_{-1}(\bar x_{-1}+z_{-1}))+E\sum_{t=0}^T\delta_{D_t}(\bar x_t+z_t).
\]
By A, B and C, the right hand side is finite on a neighborhood of $(0,\bar u)$, so $(0,\bar u)\in\inte \dom\varphi$.  \thref{ass:rec}.1 is clear since $\dom_u Ef=\U$. As to \thref{ass:rec}.2, note first that $(z,\bar u)\in\dom Ef$ implies $z_{-1}\in\dom S_{-1}$ and $z_t\in D_t$ almost surely for $t=0,\ldots,T$. By D, there is a $u\in\U$ such that $S_{-1}(z_{-1})\le u$. By the conditional Jensen's inequality, $S_{-1}(E_tz_{-1})\le E_tu$ so $E_tz_{-1}\in\dom S_{-1}$ and thus, $S_{-1}(E_tz_{-1})\in\U$, by D. Applying Jensen's inequality (see e.g.\ \cite[Theorem~8]{pp22b}) to the indicator functions of $D_t$ shows that $E_tz_{t'}\in D_{t'}$ for $t'\le t$. The process $\bar z$ in \thref{ass:rec}.2 can thus be taken $\bar z_{t'}=E_tz_{t'}$ for $t'\le t$ and $\bar z_{t'}=0$ for $t'>t$.
\end{proof}

If, in addition to \thref{ass:fmds}, we assume that there are no portfolio constraints and that $S_{-1}$ is finite on $\reals^K$, the right side of the inequality in the above proof is, by \thref{cor:IFmc}, $\tau(\X\times\U,\V\times\Y)$-continuous under \thref{ass:fmds}. Since a convex function that is bounded from above on an open set is continuous on the set, this implies Mackey continuity of $\varphi$ and thus the existence of a subgradient in $\V\times\Y$. In that case, one can thus avoid going through the arguments in Sections~\ref{sec:sdp} and \ref{sec:recourse}. In the presence of portfolio constraints, the Mackey continuity fails but we still find a dual solution in $\V\times\Y$. Note also that, the above proof gives the existence of a solution in the strong dual \eqref{sd} even without part D in \thref{ass:fmds}. Strong dual solutions would exist even if, in part C, the finiteness of $EV$ was weakened by only requiring that $EV$ be strongly continuous at $\bar u-\sum_{t=0}^{T-1}\bar x_t\Delta s_{t+1}-C\cdot \bar x_{-1}+S_{-1}(\bar x_{-1})$; see \thref{ex:orliczcont}.





\section{Appendix}\label{sec:appendix}

\subsection{Relative core of a convex set}\label{sec:hulls}

Let $U$ be a vector space and $C\subset U$ convex. The {\em positive hull}
\[
\pos C := \bigcup_{\lambda>0}\lambda C
\]
of $C$ is the intersection of all cones containing $C$. If $0\in C$ and $0<\lambda_1<\lambda_2$, then by convexity, $\lambda_1 C\subseteq\lambda_2 C$.

\begin{lemma}\thlabel{lem:hulls}
Given a linear $A:X\to U$ and convex sets $C,C'\subset X$ and $D\subset U$, we have 
\begin{enumerate}
\item
  $\pos(AC)=A\pos C$,
\item
  $\pos(A^{-1}D)=A^{-1}\pos D$,
\item
$\pos(C\times D)=\pos C\times\pos D$ if $0\in C$ and $0\in D$,
\item
  $\pos(C\cap C')=\pos C\cap\pos C'$ if $0\in C$ and $0\in C'$,
\item
  $\pos(C+C')=\pos C+\pos C'$ if $0\in C$ and $0\in C'$,
\end{enumerate}
\end{lemma}

\begin{proof}
The first two claims are clear. As to 3, since $0\in C$ and $0\in D$, we have $C\times D\subset\pos C\times\pos D$ so $\pos(C\times D)\subset\pos C\times\pos D$. If $(x,u)\in\pos C\times\pos D$, we have $x\in\lambda_1 C$ and $u\in\lambda_2 D$ for some $\lambda_i>0$. Since the sets contain the origins, we have $x\in\max\{\lambda_1,\lambda_2\}C$ and $u\in\max\{\lambda_1,\lambda_2\}D$, so $(x,u)\in\pos(C\times D)$. Defining $A:X\to X\times X$ by $Ax=(x,x)$, we have $C\cap C'=A^{-1}(C\times C')$, so 4 follows from 2 and 3. Defining $A:X\times X\to X$ by $A(x,x')=x+x'$, we have $A(C\times C')=C+C'$, so 5 follows from 1 and 3.
\end{proof}

The {\em core} of a set $C\subset U$, denoted by $\core C$, is the set of points $u\in C$ for which $\pos(C-u)= U$. The {\em relative core} of $C$, denoted by $\rcore C$, is the core of $C$ relative to the affine hull of $C$. Recall that the {\em affine hull} $\aff C$ of $C$ is the smallest affine set containing $C$. A set $C$ is {\em affine} if $\lambda u+(1-\lambda)u'\in C$ for all $u,u'\in C$ and $\lambda\in\reals$.

\begin{lemma}\thlabel{lem:rcore}
Given a convex set $C$,
\begin{enumerate}
\item
  $\rcore C=\{u\in C\mid \pos(C-u)=\aff(C-u)\}$,
\item
  $\pos(C-x)=\pos(C-x')$ for every $x,x'\in\rcore C$.
\end{enumerate}
\end{lemma}

\begin{proof}
By definition, $\pos(C-u)\subseteq\aff(C-u)$ for any $u\in C$. The converse holds if and only $u\in\rcore C$. Since $\aff(C-u)$ is independent of the choice of $u\in C$, the second claim follows from the first one.
\end{proof}

\subsection{Continuity of convex functions}

Let $U$ be a topological vector space.

\begin{theorem}\thlabel{thm:cont0}
A convex function which is bounded from above on an open set is either proper and continuous or identically $-\infty$ throughout the core of its domain. 
\end{theorem}

\begin{proof}
In the case of a proper convex function, the proof can be found e.g.\ in \cite[Proposition~I.2.5]{et76}. A simple line segment argument shows that, if a convex function equals $-\infty$ at some point, then it equals $-\infty$ throughout the core of its domain.
\end{proof}

We say that a function $g$ is {\em relatively continuous} at $u\in\dom g$ if $g$ is continuous at $u$ relative to $\aff\dom g$. The following a straightforward consequence of \thref{thm:cont0}.

\begin{corollary}\thlabel{cor:cont0}
A convex function which is bounded from above on a relatively open set of $\aff\dom g$ is either proper and relatively continuous or identically $-\infty$ throughout $\rcore\dom g$. 
\end{corollary}

Even simpler conditions for continuity are available if the underlying space is {\em barreled} in the sense that $\core C=\inte C$ for every closed convex set $C$. Fr\'echet spaces are barreled. The following implies that the same holds if, instead of being closed, $C$ is the domain of a lsc convex function.

\begin{theorem}\thlabel{lscbar}\cite[Corollary~8B]{roc74}
In a barreled space, a lsc convex function is either proper and continuous or identically $-\infty$ throughout the core of its domain.
\end{theorem}

The {\em relative interior} $\rinte C$ of a set $C$ is the interior of $C$ with respect to $\aff C$. If $C$ is a closed subset of a topological vector space $U$ and $\aff C$ is barreled, then the earlier argument applied on (the linear translation of) $\aff C$ gives
\[
\rinte C=\rcore C.
\]

The following is an immediate corollary of \thref{lscbar}. We say that an affine set is {\em barreled} if its translation to the origin is barreled. Again, a closed affine set in a metrizable space is barreled.

\begin{corollary}\thlabel{cor:lscbar}
Let $g$ be a lsc convex function on a topological vector space. If $\aff\dom g$ is barreled, then $g$ is either proper and relatively continuous or identically $-\infty$ throughout $\rcore\dom \varphi$.
\end{corollary}


\begin{remark}\thlabel{fb}
Closed subspaces of barreled spaces need not be barreled but closed subspaces of Fr\'echet spaces are Fr\'echet, and thus, barreled.
\end{remark}

\begin{corollary}\thlabel{cor:relcontiff}
A proper lsc convex function $g$ with $\rcore\dom g\ne\emptyset$ on a Fr\'echet space is relatively continuous throughout $\rcore\dom g$ if and only if $\aff\dom g$ is closed.
\end{corollary}

\begin{proof}
A closed subspace of a Fr\'echet space is Fr\'echet and, in particular, barreled. Thus, by \thref{cor:cont0}, the closedness of $\aff\dom g$ implies relative continuity. To prove the converse, let $(u^\nu)_{\nu=\infty}^\infty$ be a sequence in $\aff\dom g$ converging to a $\bar u$. The sequence is bounded so, for every neighborhood $N$ of the origin, there is an $\alpha>0$ such that $(u^\nu)_{\nu=\infty}^\infty\subset\alpha N$. Relative continuity of $g$ gives the existence of an open set $N$, a $u_0\in\rcore\dom g$ and an $\beta>g(u_0)$ such that $N\cap\aff\dom g\subset(\lev_\beta g-u_0)$. It follows that there is a constant $\alpha>0$ such that $g(\alpha u^\nu)\le\beta$ for all $\nu$. The lower semicontinuity of $g$ gives $g(\alpha\bar u)\le\beta$ so $\bar u\in\aff\dom g$.
\end{proof}

The following extends \thref{cor:lscbar} by relaxing the lower semicontinuity assumption. 

\begin{theorem}\thlabel{thm:ipcont}\cite[Theorem 2.7.1.(vi)]{zal2}
Let $X$ be a Fr\'echet space, $F$ a lsc convex function on $X\times U$ and
\[
\varphi(u):=\inf_{x\in X}F(x,u).
\]
If $\aff\dom\varphi$ is barreled, then $\varphi$  is either proper and continuous or identically $-\infty$ throughout $\rcore\dom \varphi$.
\end{theorem}

\subsection{Conjugates and subgradients}\label{sec:cs}

Let $U$ and $Y$ be vector spaces in {\em separating duality} under a bilinear form
\[
(u,y)\mapsto\langle u,y\rangle.
\]
The {\em conjugate} of $g: U\to\ereals$ is the extended real-valued convex function on $\Y$ defined by
\[
g^*(y)=\sup_{u\in U}\{\langle u,y\rangle-g(u)\}.
\]
The conjugate of the {\em indicator function}
\[
\delta_C(u):=
\begin{cases}
  0 & \text{if $u\in C$},\\
  +\infty &\text{otherwise}
\end{cases}
\]
of a set$C\subset U$ is the {\em support function}
\[
\sigma_C(y):=\sup_{u\in C}\langle u,y\rangle
\]
of $C$.

By definition, a function and its conjugate satisfy the {\em Fenchel inequality}
\[
g(u)+g^*(y)\ge\langle u,y\rangle
\]
for all $u\in U$ and $y\in Y$. When $g(u)$ is finite, the inequality holds as an equality iff
\[
g(u')\ge g(u)+\langle u'-u,y\rangle\quad\forall u'\in U.
\]
We then say that $y$ is a {\em subgradient} of $g$ at $u$. The set of subgradients of $g$ at $u$ is known as the {\em subdifferential} of $g$ at $u$ and denoted by $\partial g(u)$. The subdifferential is defined as the empty set unless $g(u)$ is finite.

\begin{theorem}\thlabel{relsd}\cite[Theorem 2.4.12]{zal2}
If $g$ is relatively continuous and finite at $u\in\dom g$, then $\partial g(u)\ne\emptyset$.
\end{theorem}

Given a proper convex function $g$, its {\em recession function} is the convex positively homogeneous function given by
\[
g^\infty(u)=\lim_{\alpha\upto\infty}\frac{g(\bar u+\alpha u)-g(\bar u)}{\alpha},
\]
where $\bar u\in\dom g$; see e.g.\ \cite{roc66}.

\begin{lemma}\thlabel{domrec}\cite[Corollary~3C]{roc66}
If $g^*$ is proper, then
\[
(g^*)^\infty=\sigma_{\dom g}.
\]
\end{lemma}

The following was used in the proof of \thref{lem:spi}.

\begin{lemma}\thlabel{lem:idemp}
Let $\pi$ be a continuous linear mapping on $U$ and $g:U\to\ereals$ such that $g\circ\pi\le g$. Then $g^*\circ\pi^*\le g^*$. If $\pi$ is idempotent and $g$ is subdifferentiable at $u\in\rge\pi$ in the relative topology of $\rge\pi$, then $g$ is subdifferentiable at $u$.
\end{lemma}

\begin{proof}
The first claim follows from
\begin{align*}
g^*(\pi^* y) &=\sup\{\langle u,\pi^*y\rangle -g(u)\}\\
&\le \sup\{\langle\pi u,y\rangle -g(\pi u)\}\\
&\le \sup\{\langle u,y\rangle -g(u)\}\\
&= g^*(y).
\end{align*}
Assume now that $g$ is subdifferentiable at $u\in\rge\pi$ relative to $\rge\pi$. By Hahn-Banach, linear functionals on $\rge\pi$ can be expressed by elements of $Y$. Thus, there is $y\in Y$ such that
\[
g(u') \ge g(\pi u')\ge g(u) + \langle \pi u'-u,y\rangle = g(u) + \langle u'-u,\pi^*y\rangle\quad\forall u'\in U,
\]
where the last equality holds when $\pi$ is idempotent.
\end{proof}

\subsection{Duality in optimization}\label{sec:cd}

This appendix gives a brief summary of the conjugate duality framework of \cite{roc74}. We include an extra linear perturbation in the primal objective which allows us to formulate certain results in a more convenient form. It does not, however, interfere with the original arguments from \cite{roc74} as we see below.

Let $X$ and $U$ be in separating duality with $V$ and $Y$, respectively. Given a convex function $F$ on $X\times U$, the primal problem is
\[
\minimize\quad F(x,u)-\langle x,v\rangle\quad\ovr x\in X.
\]
and the dual problem is
\[
\maximize\quad \langle u,y\rangle-F^*(v,y)\quad\ovr y\in Y.
\]
The primal value function is
\[
\varphi_v(u):=\inf_{x\in X} \{F(x,u) - \langle x,v\rangle\}
\]
and the dual value function
\[
\gamma_u(v):=\inf_{y\in Y} \{F^*(v,y)-\langle u,y\rangle\}.
\]
By Fenchel's inequality, 
\[
F(x,u)+F^*(v,y)\ge\langle x,v\rangle+\langle u,y\rangle
\]
for all $(x,u)\in X\times U$ and $(v,y)\in V\times Y$, so
\[
\varphi_v(u)\ge -\gamma_u(v)\quad\forall u\in U, v\in V.
\]
If $\varphi_v(u)> -\gamma_u(v)$, a {\em duality gap} is said to exist.

By definition, $\varphi_v^*(y)=F^*(v,y)$, so the dual problem can be written as
\[
\maximize\quad \langle u,y\rangle-\varphi_v^*(y)\quad\ovr y\in Y.
\]
The properties of conjugates and subgradients from the previous section thus imply that the absence of a duality gap and existence of dual solutions come down to closedness and subdifferentiability of $\varphi_v$ at $u$; see \thref{ndg,sdvarphi} below.

The {\em Lagrangian} associated with $F$ is the convex-concave function on $X\times Y$ given by
\[
L(x,y):=\inf_{u\in U}\{F(x,u)-\langle u,y\rangle\}.
\]
Clearly, the conjugate of $F$ can be expressed as
\[
F^*(v,y) = \sup_{x\in X}\{\langle x,v\rangle - L(x,y)\}.
\]
The Lagrangian saddle-point problem is to find a saddle-value and/or a saddle-point of the convex-concave function
\[
L_{v,u}(x,y) := L(x,y)-\langle x,v\rangle+\langle u,y\rangle.
\]
We have
\[
\langle u,y\rangle - F^*(v,y) = \inf_xL_{v,u}(x,y)
\]
so the dual problem is equivalent to the maximization half of the Lagrangian minimax problem. When $F$ is closed in $u$, the biconjugate theorem gives
\[
F(x,u)-\langle x,v\rangle = \sup_yL_{v,u}(x,y)
\]
so the primal problem is the minimization half of the minimax problem.

\begin{theorem}\thlabel{ndg}
The following are equivalent,
\begin{enumerate}
\item There is no duality gap. 
\item $\varphi_v$ is closed at $u$. 
\end{enumerate}
If $F$ is closed in $u$, the above are equivalent to
\begin{enumerate}[resume]
\item $L_{v,u}$ has a saddle-value.
\end{enumerate}
If $F$ is closed, the above are equivalent to
\begin{enumerate}[resume]
\item $\gamma_u$ is closed at $v$. 
\end{enumerate}
\end{theorem}

\begin{proof}
We have $\varphi_v^*(y)=\sup_{v,y}\{\langle x,v\rangle+\langle u,y\rangle -F(x,u)\}=F^*(v,y)$, so
\[
\varphi_v^{**}(u)=\sup_y \{\langle u,y\rangle-F^*(v,y)\} =-\gamma_u(v),
\]
and the equivalence of 1 and 2 follows from the biconjugate theorem. When $F$ is closed in $u$, 1 is equivalent to 3 by the remarks before the statement.
When $F$ is closed, $F=F^{**}$ by the biconjugate theorem, so 4 is equivalent to 1, by symmetry.
\end{proof}

\begin{theorem}\thlabel{sdvarphi}
If $\varphi_v(u)<\infty$, then the following are equivalent,
\begin{enumerate}
\item
  There is no duality gap and $y$ solves the dual.
\item
  Either $\varphi_v(u)=-\infty$ or $y\in\partial\varphi_v(u)$.
\end{enumerate}
If $F$ is closed in $u$, the above are equivalent to
\begin{enumerate}[resume]
\item
  $\underset{x}\inf\,\underset{y}\sup\, L_{v,u}(x,y) = \underset{x}\inf\, L_{v,u}(x,y)$.
\end{enumerate}
If $F$ is closed, the above are equivalent to
\begin{enumerate}[resume]
\item
  $\gamma_u$ is closed at $v$ and $y$ solves the dual.
\end{enumerate}
\end{theorem}

\begin{proof}
Condition 2 means that either $\varphi_v(u)=-\infty$ or $\varphi_v(u)+\varphi_v^*(y)=\langle u,y\rangle$. The equivalence of 1 and 2 thus follows from the biconjugate theorem. When $F$ is closed in $u$, the equivalence of 1 and 3 follows from the equivalence of 1 and 3 in \thref{ndg}. When $F$ is closed, the equivalence of 1 and 4 follows from the equivalence of 1 and 4 in \thref{ndg}. 
\end{proof}

\begin{theorem}\thlabel{lem:sgvarphi}
The following are equivalent,
\begin{enumerate}
\item There is no duality gap, $x$ solves the primal, $y$ solves the dual and both problems are feasible.
\item $y\in\partial\varphi_v(u)$ and $x$ solves the primal.
\item $(v,y)\in\partial F(x,u)$.
\end{enumerate}
If $F$ is closed in $u$, the above are equivalent to
\begin{enumerate}[resume]
\item $v\in\partial_xL(x,y)$ and $u\in\partial_y[-L](x,y)$.
\end{enumerate}
If $F$ is closed, the above are equivalent to 
\begin{enumerate}[resume]
\item $x\in\partial\gamma_u(v)$ and $y$ solves the dual.
\item $(x,u)\in\partial F^*(v,y)$.
\end{enumerate}
\end{theorem}

\begin{proof}
The equivalence of 1 and 2 follows from \thref{sdvarphi}. Condition 2 means that
\[
F(x',u') - \langle x',v\rangle \ge F(x,u)- \langle x,v\rangle+\langle u'- u,y\rangle \quad\forall x',u'
\]
which is 3. It is clear from the discussion just before \thref{ndg} that, when $F$ is closed in $u$, 1 means that $(x,y)$ is a saddle-point of $L_{v,u}$, which is 4. 
When $F$ is closed, 3 is equivalent to 6, by the biconjugate theorem. The equivalence of 5 and 6 follows like that of 2 and 3.
\end{proof}

\begin{theorem}\thlabel{lem:cdrec}
Assume that $v\in V$ is such that $\varphi_v$ is either relatively continuous at $u\in U$ or $\varphi_v(u)=-\infty$. Then
\[
\gamma_u(v):=\inf_{y\in Y} \{F^*(v,y)-\langle u,y\rangle\}
\]
is closed at $v$ and
\[
\inf_x \{F(x,u)-\langle  x,v \rangle \}= \sup_y\{\langle u,y \rangle - F^*(v,y)\}, 
\]
where the supremum is attained. 
If $F$ is proper lsc and the first assumption holds for all $v\in V$, then $\gamma_u$ is proper lsc, and
\begin{equation}\label{eq:gammarec}
\gamma_u^\infty(v)=\inf_{y\in Y}\{(F^*)^\infty(v,y)-\langle u,y\rangle\},
\end{equation}
where the infimum is attained.
\end{theorem}

\begin{proof}
If $\varphi_v(u)$ is finite then, by \thref{relsd}, there exists $y\in\partial\varphi_v(u)$ which gives the second claim. The claims in the first paragraph follow now from  \thref{sdvarphi}.

Under the additional assumption, $\gamma_u$ is proper lsc, so, by the biconjugate theorem, $F(\cdot,u)^*=\gamma_u$. Since $F$ is proper lsc, there exists $(\bar v,\bar y)\in\dom F^*$. By Fenchel,
\[
F(x,u)+F^*(\bar v,\bar y)\ge \tilde F(x,u),
\]
where $\tilde F(x,u) := \delta_{\dom F}(x,u)+\langle x,\bar v\rangle + \langle u,\bar y\rangle$. Thus,
\[
\varphi_v(u)+F^*(\bar v,\bar y) \ge \tilde\varphi_v(u),
\]
where
\[
\tilde\varphi_v(u):=\inf_x\{\tilde F(x,u)-\langle x,v\rangle\}.
\]
We have $\dom\tilde\varphi_v=\dom\varphi_v$ and, by assumption, $u\in\dom\varphi_v$. When $\tilde \varphi_v(u)$ is finite, the above inequality implies that $\varphi_v(u)$ is finite as well, and then continuity of $\varphi_v$ at $u$ implies that of $\tilde\varphi_v$. Thus, by the first part,
\[
\sup_x \{\langle x,v\rangle-\tilde F(x,u)\} = \inf_{y\in Y}\{\tilde F^*(v,y)-\langle u,y\rangle\},
\]
where the infimum is attained. Since $F(\cdot,u)^*=\gamma_u$, \thref{domrec} implies $\tilde F(\cdot,u)^*(v)=\gamma_u^\infty(v-\bar v)-\langle u,\bar y\rangle$ and $\tilde F^*(v,y)= (F^*)^\infty(v-\bar v,y-\bar y) $. Substituting into above,
\begin{align*}
\gamma_u^\infty(v-\bar v)-\langle u,\bar y\rangle &= \inf_{y\in Y}\{(F^*)^\infty(v-\bar v,y-\bar y)-\langle u,y\rangle\}\\
&  = -\langle u,\bar y\rangle+ \inf_{y\in Y}\{(F^*)^\infty(v-\bar v,y)-\langle u,y\rangle\},
\end{align*}
which proves the claim.
\end{proof}

Note that $\dom\varphi_v$ does not depend on $v$. We denote this set by $\dom\varphi$. \thref{thm:ipcont} gives the following.

\begin{theorem}\thlabel{thm:core}
If $X$ and $U$ are Fr\'echet, $F$ is lsc and $\aff\dom\varphi$ is closed, then for every $v\in V$, $\varphi_v$ is either proper and continuous or identically $-\infty$ throughout $\rcore\dom \varphi$.
\end{theorem}

The conditions of \thref{thm:core} clearly hold if $0\in\core\dom\varphi$. This, in turn, generalizes the classical Slater-type conditions which require $0\in\inte\dom\varphi$. The primal problem is said to be {\em strictly feasible} if there is an $\bar x\in X$ such that $\pos(\dom F-(\bar x,\bar u))$ is linear and closed. Since $\dom\varphi$ is the projection of $\dom F$ to $U$, \thref{lem:hulls} gives
\begin{align*}
\pos(\dom\varphi-\bar u) &= \{u\mid \exists x:\ (x,u)\in\pos(\dom F-(\bar x,\bar u))\}.
\end{align*}
Thus, the assumptions of \thref{thm:core} are satisfied if the primal is strictly feasible and the projection of the closed linear subspace $\aff\dom F$ is closed in $U$. The closedness is a classical question in functional analysis and various conditions have been given. We refer the reader to Section~\ref{sec:sdp} for examples in the context of stochastic optimization.

Combining \thref{lem:cdrec,thm:core} gives the following.

\begin{theorem}\thlabel{domcd}
If $X$ and $U$ are Fr\'echet, $F$ is lsc and $\aff\dom\varphi$ is closed, then for every $v\in V$ and $u\in\rcore\dom\varphi$, the infimum in the definition of $\gamma_u$ is attained, $\gamma_u$ is lsc proper and
\begin{align*}
\gamma_u^\infty(v)=\inf_{y\in Y}\{(F^*)^\infty(v,y)-\langle u,y\rangle\},
\end{align*}
where the infimum is attained.
\end{theorem}

\begin{proof}
By \thref{thm:core}, $\varphi_v$ is relatively continuous for every $v\in V$. The claims thus follow from \thref{lem:cdrec} and recalling that, by \thref{sdvarphi}, the subgradients of $\varphi_v$ at $u$ are exactly the minimizers in the definition of $\gamma_u$.
\end{proof}

The following was proved for Banach spaces in \cite{ab86}.

\begin{corollary}[Attouch-Brezis]\thlabel{cor:infcomp}
Let $f_1,f_2$ be closed convex functions on a Fr\'echet space $X$ such that $\aff(\dom f_1-\dom f_2)$ is closed and
\[
0\in\rcore(\dom f_1-\dom f_2).
\]
Then
\[
(f_1+f_2)^*(v)=\inf_{y\in X^*}\{f_1^*(v-y)+ f_2^*(y)\}\quad\forall v\in X^*,
\]
where the infimum is attained.
\end{corollary}

\begin{proof}
We apply the general duality framework with $F(x,u)=f_1(x)+f_2(x+u)$. We have
\[
\dom\varphi_v = \{u\mid \exists x\in\dom f_1:\ u+x\in\dom f_2\} = \dom f_2-\dom f_1
\]
and
\begin{align*}
F^*(v,y) &= \sup_{x,u}\{\langle x,v\rangle + \langle u,y\rangle - f_1(x) - f_2(x+u)\}\\
&= \sup_{x,u'}\{\langle x,v\rangle + \langle u'-x,y\rangle - f_1(x) - f_2(u')\}\\
&=f_1^*(v-y)+f_2^*(y)
\end{align*}
so the claims follow from \thref{domcd}.
\end{proof}

The following was used in \thref{ex:scaled}.

\begin{example}\thlabel{finalU}
Let $A$ be a linear mapping from a Fr\'echet $X$ space to a vector space $U$ such that $\ker A$ is closed. Then $\rge A$ is a Fr\'echet space with respect to the final topology induced by $A$. Moreover, the topological dual of $\rge A$ can be identified with $(\ker A)^\perp$ in the sense that for every $u^*\in(\rge A)^*$ there is a unique $v\in(\ker A)^\perp$ such that
\[
\langle Ax,u^*\rangle=\langle x,v\rangle
\]
for all $x\in X$.
\end{example}

\begin{proof}
The quotient space $X/\ker A$ and $\rge A$ are in one-to-one correspondence under the linear mapping $i$ that sends an $[x]\in X/\ker A$ to $Ax$. Indeed, $i^{-1}(u)=\{x\in X\mid Ax=u\}$. The final topology on $\rge A$ is topologically isomorphic under $i$ with the quotient space topology on $X/\ker A$. If $\ker A$ is closed, then $X/\ker A$ is Hausdorff and $X/\ker A$ is Fr\'echet and so too is $\rge A$.

Continuous linear functionals $l$ on $\rge A$ can be expressed as $l(u)=l'(i^{-1}(u))$, where $l'$ is a continuous linear functional on $X/\ker A$. The continuous linear functionals $l'\in(X/\ker A)^*$ are given by $l'([x])=\langle x,v\rangle$, where $v\in(\ker A)^\perp$.
\end{proof}

The following was used in the proof of \thref{lem:pircore}.

\begin{lemma}\thlabel{lem:sss}
Let $L\subset X$ and $A: X\to U$ a continuous linear mapping such that $AL$ is closed. Then $\ker A+L$ is closed. In particular, if $L$ is closed and linear and $\pi$ a continuous linear idempotent mapping on $U$ such that $\pi L\subseteq L$, then $\rge \pi +L$ is closed.
\end{lemma}
\begin{proof}
We have
\begin{align*}
\ker A+L &=\{z+x'\mid Az=0,\, x'\in L\}\\
&=\{x\mid A(x-x')=0,\, x'\in L\}\\
& = \{x\in X\mid Ax\in AL\}\\
& = A^{-1}(AL),
\end{align*}
which proves the first claim. The second claim follows from the first by choosing $A=I-\pi$. Indeed, $AL= L-\pi L=L$ is closed and, since $\pi$ is idempotent, $\ker A=\rge\pi$.
\end{proof}

\bibliographystyle{plain}
\bibliography{sp}

\end{document}